\documentclass[11pt,oneside,fleqn,reqno,titlepage]{article}
\usepackage{amsmath}
\usepackage{amssymb}
\usepackage{amsthm}
\usepackage{natbib} \citeindextrue
\usepackage{comment}
\usepackage{url}
\usepackage{bbm}
\usepackage{bm}
\usepackage{multirow}
\usepackage{relsize}
\usepackage[pdftex]{graphicx}
\usepackage{color}

\usepackage{mathrsfs}
\usepackage{amsfonts}
\usepackage{dsfont}
\usepackage{algorithmic}
\usepackage{algorithm}

\usepackage{latexsym}
\usepackage{float}
\usepackage{flushend}
\usepackage{cuted}
\usepackage{color}

\usepackage{endnotes}
\let\footnote=\endnote

\usepackage{natbib}
 \bibpunct[, ]{(}{)}{,}{a}{}{,}%
 \def\bibsep{\smallskipamount}%
 \def\newblock{\ }%

\newcommand{\argmin}{{\rm arg}\min}

\newcommand{\doublespace}{\addtolength{\baselineskip}{.25\baselineskip}}

\newcommand{\singlespace}{\addtolength{\baselineskip}{-.5\baselineskip}}

\begin{document}


\title{The Information-Collecting Vehicle Routing Problem: Stochastic Optimization for Emergency Storm Response}
\author{Lina Al-Kanj, Warren B. Powell and Belgacem Bouzaiene-Ayari}
\date{\today}

\maketitle

\clearpage

\begin{abstract}
{\color{black}
We address the problem of mitigating damage to a power grid following a storm by managing a vehicle that has to be routed while simultaneously performing two tasks: learning about damage from the grid (which requires direct observation) and repairing damage that it observes.  The learning process is assisted by calls from customers notifying the utility that they have lost power (``lights-out calls'').  However, when a tree falls and damages a line, it triggers the first upstream circuit breaker, which results in power outages for everyone on the grid below the circuit breaker.  We present a dynamic routing model that captures observable state variables such as the location of the truck and the state of the grid on segments the truck has visited, and beliefs about outages on segments that have not been visited.  Trucks are routed over a physical transportation network, but the pattern of outages is governed by the structure of the power grid.  We introduce a form of Monte Carlo tree search based on information relaxation that we call {\it optimistic MCTS} which improves its application to problems with larger action spaces.  We show that the method significantly outperforms standard escalation heuristics used in industry.}
\end{abstract}

\clearpage
\pagestyle{empty}
{\singlespace
\tableofcontents
}
\clearpage

\setcounter{page}{1}
\pagestyle{plain}

\doublespace

\section{Introduction}
Climate change is producing more powerful storms, increasing the frequency and severity of outages in the power grid (\citet{Konisky2016}). On average $55$\% of power outages in the U.S. are due to weather and it can reach up to 80\% in some years~(\citet{Camp12}).  Wind and ice bring trees and branches down on power lines creating sporadic outages that quickly spread through the grid due to the limited number of protective devices that are triggered due to a short circuit. Despite the importance of electricity in every aspect of our lives, we often do not know the location of  the fault causing an outage, which complicates the task of restoring power. Instead, utilities depend heavily on phone calls (known as ``lights-out calls'') from customers who have lost their power, in addition to information from the tripping of some circuit breakers/protective devices; however, only a few percent of customers call when they lose their power, creating tremendous uncertainty in the knowledge of the state of the grid. This in turn complicates the task of dispatching utility trucks to the location of faults (which are unknown to the utility center) to restore the grid as quickly as~possible.

\begin{figure}[tb]
\begin{center}
    \includegraphics[width=4.75in]{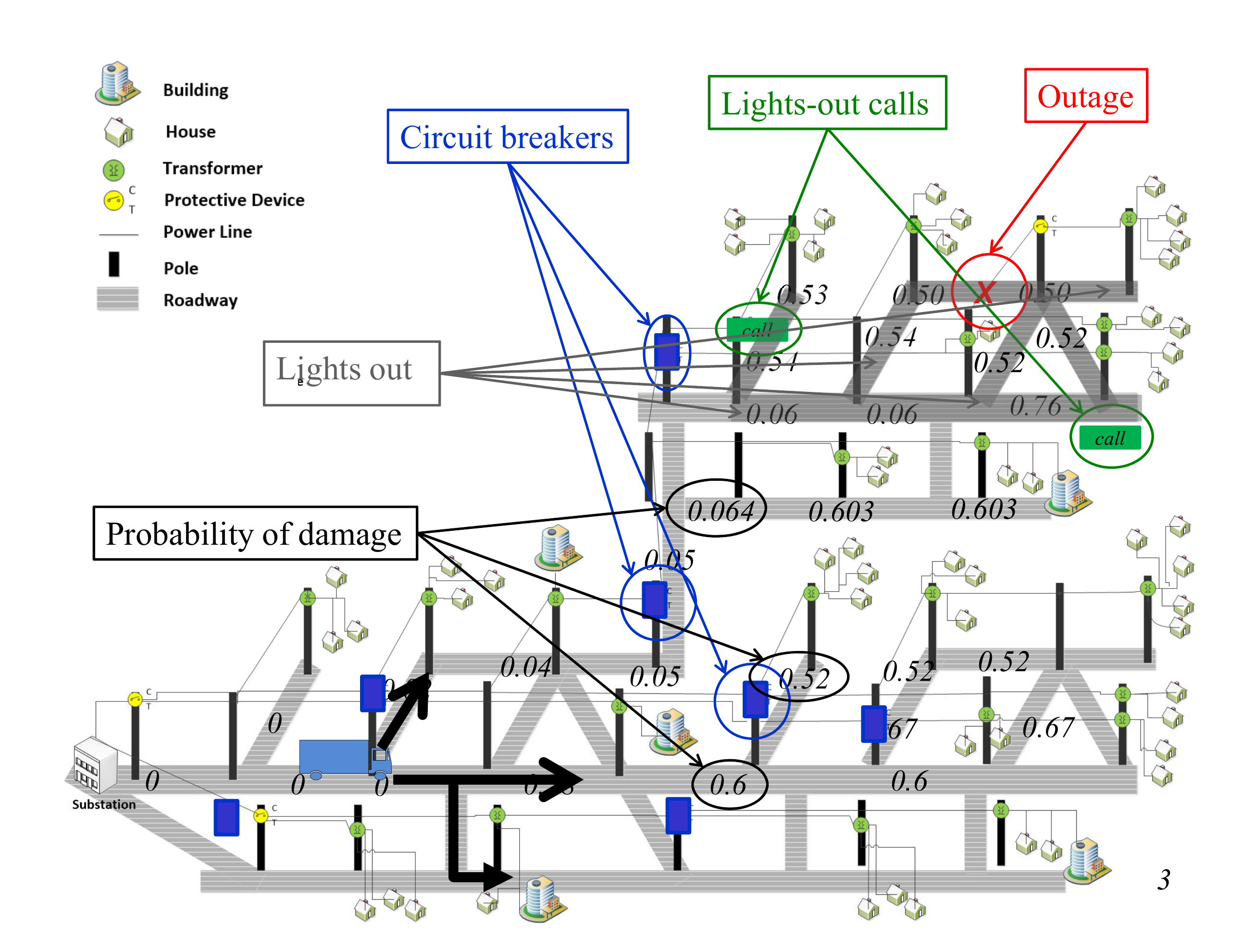}
    \caption{Illustration of grid showing lights-out calls, an outage, circuit breakers, probability of where outages might have occurred, and the segments of the grid where customers have lost their power.}
    \label{fig:stormresponseillustration}
\end{center}
\end{figure}

{\color{black}
Figure \ref{fig:stormresponseillustration} illustrates the situation.  Imagine that a tree has fallen down on a line creating an outage (we show one, but there could be several), but we do not know where it is.  The outage blows the first upstream circuit breaker, which results in everyone downstream losing their power.  Some percentage (possibly as low as one percent) will call in reporting their lights are out, but these calls typically do not come from the grid segment where the outage has occurred.  As the truck moves, it determines whether an outage has occurred on the grid segments it observes from the road network; if it observes an outage, it stops and performs a repair, which is random because it depends on the type of repair required.  It will then use this information, along with new lights-out calls, to update its belief about where outages might be.  We assume the truck knows the structure of the grid and the location of circuit breakers and this is used to update beliefs.  The challenge is to decide how to route the truck, where it has to strike a balance between moving along a segment where an outage might have occurred, and moving along a segment just to learn whether an outage has occurred.
}

Even if the EUC knows the locations of the faults, routing the truck to restore the grid to minimize the number of customers in outage is a challenge.  Currently, EUCs use a simple policy to route the utility trucks to restore power, such as first-call, first-served, which is far from optimal. In this paper, a utility truck is routed across the grid while collecting information on its way. This information is then used to update the belief model about the locations of faults, where ``belief'' refers to an estimate of the probability that an outage has occurred to the grid along a link of the network.  Now, we have a physical state (location of the truck), information state (outage calls) and belief state (the probability that a segment of the grid is out).  {\color{black} We propose a stochastic lookahead policy, and use Monte Carlo tree search (MCTS) to produce a policy that is an asymptotically optimal solution of the lookahead.

We show how to formulate our problem as a sequential decision problem which requires a high dimensional state variable, a property that raises questions about the ``curse of dimensionality.'' MCTS does not require enumerating the full state space.  Rather, it explores the future using a combination of a lookahead policy and Monte Carlo sampling.  The method is sensitive to the number of decisions that can be made at each point in time, which is a reason why we are limited to a single truck.  However, the truck does not decide whether to move left, right or straight on the road network; rather, it decides which location on the grid to investigate next, which is much larger.  We use a new version of MCTS which we call {\it optimistic MCTS} which explores new nodes in the tree using an estimate that is biased downward (we are minimizing costs) which we have found works better with large action spaces than the conventional form of MCTS (widely used in computer science) that uses a pessimistic estimate.

This paper focuses on the problem of storm response, an issue of increasing importance with the more powerful storms produced by climate change.  However, there are many application areas that involve the management of discrete resources in the presence of online learning.  Some examples include:
\begin{itemize}
  \item Disease mitigation - We may manage a resource (e.g. a medical team) that moves around a region, where it is necessary to determine a) when to stop and perform testing and b) when to offer treatments.
  \item Invasive plant species - A botanist might move around a field looking for evidence of an invasive plant species and taking steps to reduce the species (e.g. herbicides).
  \item Fire fighting - A helicopter might fly around a region looking for fire, and then dropping retardants as needed.
  \item Ocean oil contamination - A ship may travel around a body of water (e.g. Gulf of Mexico, the Mediterranean) looking for oil contamination, which requires determining where to stop and take samples, and when to remove surface pollution.
\end{itemize}
The modeling and algorithmic framework presented in this paper can be applied to any of these problems
}

This paper makes the following contributions:
\begin{itemize}
\item[1)] It provides a formal model of a stochastic vehicle routing problem that explicitly captures the physical state (the location of the truck and known state of the grid), information state (e.g. outage calls), and belief state (the probabilities that segments of the grid are out).
\item[2)] {\color{black} It proposes a stochastic lookahead policy that is solved using a new implementation of Monte Carlo tree search we call {\it optimistic} Monte Carlo tree search.  It offers an asymptotically optimal solution to the stochastic lookahead model, which is a first for the field of stochastic, dynamic vehicle routing.}
\item[3)] Our implementation of MCTS is novel, in that it replaces the use of a suboptimal heuristic rollout policy for approximating the value of new nodes in the lookahead tree with an optimistic estimate obtained by optimally solving a sampled approximation, which is a form of information relaxation.  This work motivated the convergence proof which is now reported in \cite{JiangMCTS2020}.
\item[4)] The performance of our stochastic lookahead policy is tested using a simulator that models the power grid of the state of New Jersey by using real data provided by PSE\&G which is the largest utility in New Jersey. Simulations show that a) the lookahead learning policy performs well against the posterior optimal, and b) significantly outperforms standard industry heuristics when tested on realistic stochastic models. We also demonstrate that the industry heuristic does not make good use of better information, while our method provides information-consistent behavior, with better results as information improves (as one would expect).
\end{itemize}

As a word of caution, MCTS (especially when it involves the updating of beliefs) can be computationally expensive.  There are different choices that can be made in how MCTS is implemented.  In addition, there are opportunities for sophisticated computational algorithms such as parallel processing which we did not pursue.

The paper is organized as follows. Section~\ref{sec:lit_rev} summarizes the literature of fault identification, utility crew dispatching for restoring power distribution systems, and stochastic vehicle routing. Section~\ref{sec:sys_model} gives a mathematical model of the flow of information for a grid as some event such as a storm evolves producing faults and loss of power. Section~\ref{sec:opt_prob} derives the sequential stochastic optimization problem that routes the utility truck across the grid and discusses several dispatch policies. Section ~\ref{sec:designingpolicies} describes the four classes of policies we can draw from, and uses this framework to outline how we have designed our policy, as well as choices used elsewhere in the literature.  Section~\ref{sec:MCTS} introduces MCTS as the lookahead policy that approximates the optimal one with good computational efficiency. Section~\ref{sec:results} compares the proposed learning policy against an industry-standard escalation policy. 

\section{Literature Review}\label{sec:lit_rev}
We first review the literature from the power community on outage prediction and the dispatching of utility crews, followed by a review of the relevant articles from the mainstream literature on dynamic vehicle routing. We close with an overview of the literature on Monte Carlo tree search.

\subsection{Fault Prediction and Utility Crew Dispatching Literature}
Utilities face two problems in responding to storm damage: 1) identifying the location of outages and 2) dispatching the utility crews.

The power grid contains protective devices that detect power flow interruption upon which they shut down power flow to the affected components to avoid further damage. Thus, in case one of the components of the circuit faults, the closest upstream protective device shuts down causing power outage to all the downstream components. Each utility designs its own approach for power restoration during and after a storm. While a few utilities have implemented advanced information systems \citep{Lamp02}, most utilities (including the focus of this paper) depend primarily on phone calls from customers who have lost their power.  In addition, only a few percent of customers make these calls, resulting in a high level of uncertainty about the location of where damage has occurred.


It is estimated that $90$\% of customer outage-minutes are due to faults affecting the local distribution systems~\citep{Camp12}. Most distribution systems use a radial structure, which means there is a single path from the substation that brings power from the transmission grid down to customers.  This structure makes it possible for utilities to design simple rules, known as escalation algorithms, which help to isolate the location of a fault~(\citet{Scot90,HLC91}).
Escalation algorithms suit single fault scenarios but cannot capture the case of multiple fault scenarios. \citet{LS99} present an improved escalation algorithm for a heat storm where the escalation from the locations of calls depends on the type of upstream devices.  However, all of these escalation algorithms are heuristics without any formal optimization theory to support them.

Artificial intelligence techniques that make use of customer calls to identify the fault locations are also investigated to provide better performance than the simple heuristics described above. For example, \citet{LTH94} present a neural network,  but are limited by the need for a large sample training set. \citet{CW98} propose an approach based on fuzzy set theory and tabu search, but the limitation is in the computational complexity that becomes intractable for large networks.  \citet{LS02} design a knowledge-based outage identification that make use of SCADA and automated meter reading to provide the EUC with knowledge about the status of the distribution system on top of customer calls. However, most utilities do not have automated distribution systems and instead rely primarily on the grid topology, the phone calls of the customers, and the experience of the utility personnel to estimate the locations of outages. Moreover, since only a few customers will call, identifying the locations of the faults across the power system is still a major dilemma for the EUC which needs to restore the network as fast as possible.  This is a challenge that all methods have to address (including ours).

Managing utility crews to restore outages attracted a modest level of attention in the research literature~(\citet{SLK16,ZDT98,Wh14}). \citet{ZDT98} routes a utility truck  to the location of a customer that called to report a power outage; then, after restoring the fault, the truck is routed to the next calling customer that has the shortest travel time. \citet{Wh14} develops a model based on data mining and machine learning techniques to predict outages in the grid using collected data from past six storms  as well as asset information (framing, pole age, etc.), in addition to environmental information. Then, given the predicted outages, deterministic optimization models are developed to route a given number of utility trucks to perform a predefined number of repair jobs required at each damaged location in order to minimize the grid restoration time. \citet{SLK16} proposes an outage prediction model at the level of areas served by a substation, and then introduces a mechanism to plan hourly crew staffing levels across different organizations (service centers, local contractors, mutual aid crews) and different crew types in order to minimize the overall grid restoration time but they do not solve the problem of utility crew routing. 

\citet{CHB11} considers the power system restoration planning problem for disaster recovery. They address the problem of routing utility trucks to restore the grid; however, the faults across the grids are modeled via scenarios and then the problem is solved via typical stochastic programming. A similar approach is applied in~\citet{HBC10} for  single commodity allocation and minimizing the delivery time of the routed vehicles post a disaster such as hurricane. These works overcome the probabilistic nature of the problem via scenario generation, but are limited by the sample chosen, and the inability to handle the fully sequential nature of the problem.  These papers are effectively introducing a static policy using an approximate stochastic lookahead (in particular, a two-stage lookahead), without modeling the adaptive process of routing and learning.  By contrast, our policy is a full multistage, stochastic lookahead policy.

Including a belief model about the location of faults in a system of interest has been addressed in the literature in a limited way. \citet{LAD16} addresses the problem of planning and scheduling maintenance operations for a set of geographically distributed machines, subject to random failures, using a two-step, iterative approach. In the first step, a maintenance model determines the optimal time until the next preventive maintenance operation. In the second step, a routing model assigns and schedules maintenance operations to each technician over the planning horizon within the workday. The problem is modelled as a stochastic combinatorial optimization problem with probabilistic constraints. \citet{FAD16} applies the same methodology but in the context of planning and scheduling preventive maintenance of sediment-related sewer blockages in a set of geographically distributed sites that are subject to nondeterministic failures.  \citet{XBN17} assumes independence of the information collection component in their problem {\color{black}(as do we)} to reduce complexity and proposes an approximate formulation to obtain an implementable policy. \citet{CEK15} addresses the stochastic debris clearance problem, which captures post-disaster situations where the limited information on the debris amounts along the roads is updated as clearance activities proceed. The problem is solved using a partially observable Markov decision process model that becomes computationally intractable even for small networks. 

This paper builds on the work in ~\citet{ABP14} which develops a Bayesian model of the probabilities that there are outages (usually from trees falling on lines) at different points in the grid based on phone calls and observations of outages.  This work uses the known structure of the grid to develop posterior beliefs from a prior and new information.

\subsection{Stochastic Vehicle Routing Literature}
The study of stochastic and dynamic vehicle routing has a long history, with survey articles appearing as early as \citet{StGo1982}.  The problem area has attracted so much attention that there has been a steady flow of reviews and survey articles (\citet{StGo1983, Ps1988, PoJa1995, GePo1998, BeCo2010, PiGe2013, LaMa2014}).  Most of this literature focuses on uncertainty in pickups or deliveries, although some authors address random travel times (\citet{KeMo2003}).   Algorithmic strategies range from scenario-based lookahead policies (e.g. (\citet{LaHa2002,BeVa2004}), to a wide variety of rolling horizon heuristics that involve solving deterministic approximations modified to handle uncertainty.

Most of the academic literature has focused on solving stochastic lookahead models (\citet{LaHa2002, BeVa2004}).  Since these problems combine uncertainty with the complexity of these difficult integer programs, these lookahead models are themselves quite hard to solve.  Often overlooked is that they are simply rolling horizon heuristics to solve a fully-sequential problem, which is typically not modeled (but is often represented in a simulator). \citet{PowellEJTL2012} provides a general framework for modeling stochastic, dynamic problems in transportation and logistics, and discussions how four classes of policies for making decisions in this context.  The proper modeling of sequential decisions in vehicle routing in terms of searching over policies is relatively new to the transportation literature, but a nice example is found in \citet{GoOh2013} which then proposes a broad class of practical roll-out policies.

{\color{black}
There are a few papers that deal with learning the state of the environment using a vehicle or drone.  \citet{Dolinskaya2018} consider an ``adaptive orienteering'' problem that involves finding a path over a stochastic network.  The optimization formulation (see equations (1) and (2) in \citet{Dolinskaya2018}) assume you choose a) which nodes to service (and the order to serve them) and then b) the path to use traversing the node, after which you learn (perfectly) the link costs.  They do not model a process of acquiring information and maintaining probabilistic beliefs that evolve as information arrives (this is central to our problem).  \cite{Maya2016} consider a network repair crew routing and scheduling problem for emergency relief distribution, which is modeled and solved as a deterministic dynamic programming problem.

\cite{Glock2020} describe the problem of a drone collecting information after an emergency (they use the example of a fire, but it could be a storm).  Although they present a nice model capturing the belief about the state of the environment, their optimization model is deterministic (see in particular their equation (16)) and they use classical deterministic search methods to optimizing the routing of the drone, which means they do not face the problem of designing a {\it policy} that makes decisions as information becomes available.

Our paper models the routing of the utility truck as a fully sequential problem that updates beliefs about the location of outages as new information arrives.  The observations the truck makes as it moves to the next location are unknown (and therefore random) when we make the decision of where to move.  Most important, our policy considers both the value of the information that we may acquire by moving the utility truck along a route, as well as the benefits from performing repairs at the same time.
}


\subsection{Monte Carlo Tree Search}
Monte Carlo tree search (MCTS) is a method for approximately solving stochastic problems over some horizon as a form of rolling horizon procedure.  The concept of MCTS traces its roots to the seminal paper~\citet{CFH05}, which describes a sampled method for solving classical Markov decision processes.  The term ``Monte Carlo tree search'' was a term coined by \cite{coulom2006}. MCTS has been popular in the computer science community (see \cite{Browne2012} and  \cite{Fu2017a} for nice reviews), but has only recently seen applications in transportation and logistics (\citet{MaNe2015, EdGa2016}), but without the formal model that we present here.

MCTS has evolved in the literature primarily in the context of deterministic problems. Deterministic MCTS has been used extensively in the literature, as summarized in~\citet{BPW12}. Stochastic outcomes (which characterizes our application) have been handled using a process known in the computer science community as ``determinization''~(\citet{BFT09,BSC07,Caz06}), whereby we explicitly represent exogenous information in the tree directly~(\citet{CHS11}), which is the approach that we take.

We describe our implementation of MCTS, which is new, in section \ref{sec:MCTS}.

{\color{black}
\subsection{Comparison to the vehicle routing literature}
This paper introduces a fundamentally new class of vehicle routing problem (information collecting VRPs), in a new application context (repairing power grids), which is solved using a new class of methods for logistics (Monte Carlo tree search), using a new class of MCTS strategy (which we call ``optimistic'' MCTS).  For this reason, we pause to put this work in the context of the transportation literature.

In the vehicle routing literature, our problem is closer to a dynamic traveling repairman problem (\cite{Garcia2002}, \cite{Luo2014}), generalized to include active learning, which is a major extension. However, we do not have to visit a specified set of cities, nor do we have to find a path from an origin to a destination. We have to move around a grid until we are sure (within some tolerance) that we do not have any outages, in a way that minimizes the total time that customers have lost their power (this objective is unique to this problem setting).  Most important is that we are not aware of any work in this area that involves active learning, which requires that the model maintain beliefs about uncertain parameters that evolve as new information arrives.

Our problem setting involves emergency response to repair grids after a storm.  \cite{CEK15} considers a debris clearance problem using a POMDP formulation (that does not scale), and proposes a heuristic lookahead policy that does not explicitly capture the value of information.  Our setting explicitly considers the setting of restoring power grids, which introduces special structure in terms of relating information on power outages to potential locations where the grid may have been damaged.  The methods are all tested on a real dataset from the largest utility in New Jersey. We also note that the problem is quite timely given the likelihood of more serious storms.

We present a formal model of the optimization problem which involves finding the best {\it policy} for dispatching a truck. There are many papers that suggest various lookahead policies (see \cite{pillac2013} for a review), without actually stating the real objective of searching over the space of policies (\cite{Goodson2013} is a notable exception, but also see \cite{PowellEJTL2012}).

From this rigorous formulation, we present for the first time a fully sequential, stochastic lookahead policy based on Monte Carlo tree search.  While MCTS is quite new to the vehicle routing community, we have introduced a novel form of MCTS that we call ``optimistic MCTS'' which solves the lookahead problem within MCTS (not to be confused with the lookahead policy we are solving with MCTS) as a deterministic optimization problem based on a sampled relaxation.  We anticipate there will be considerable experimentation with this strategy for different stochastic, dynamic vehicle routing problems.

We believe our paper is presenting the first formal model of an information-collecting VRP. It requires a state variable that includes both the physical state (the location of the vehicle), other information (e.g. phone calls), and a high-dimensional belief state giving the probability of faults on every link of the network.  A stochastic lookahead model is needed to create incentives for the vehicle to explore sections of the grid to see if outages {\it might} exist.  It is significant that MCTS is insensitive to the dimensionality of the state variable, although it is quite sensitive to the number of decisions per node, which is why we are only considering the single vehicle case in this paper.  However, while our application setting is special, the solution approach is quite general.
}

\section{Problem Description}~\label{sec:sys_model}
\noindent This section describes the dynamic model, which requires capturing two networks: the transportation network over which the utility truck is moving, and the power grid that the truck is trying to repair.  We learn indirectly about outages to the grid through sporadic calls that a customer has lost power, where we use the structure of the grid to make inferences about the location of outages. The focus of this paper is on the problem of optimizing the routing of the truck.

\subsection{The Transportation Network}
Let $G(\mathcal{V},\mathcal{E})$ be the graph representing the road network through which the utility truck can travel to check the power grid where $\mathcal{V}=\{i: i=1,\ldots,N\}$ is the set of nodes of the road network; each node in each node $\mathcal{V}$ represents either a begin/end of a road segment or a pole of the power grid. In this work, we assume an overhead distribution power grid that is mounted on poles where each pole  is represented as a node in the transportation network. Let $\mathcal{I}=\{i: i=1,\ldots,I\}$ be the set of poles  of the power grid; thus, we have $\mathcal{I}\subseteq \mathcal{V}$. Now, we define  $\mathcal{E}$ as the set containing the arcs/roads in the graph, i.e.,  $\mathcal{E}=\{(i,j): \mbox{~if there is an arc between nodes~} i \mathrm{~and~ } j \mathrm{~in~} \mathcal{V}\}$. Thus, some of the arcs/roads in $\mathcal{E}$ are parallel to the power lines of the power~grid as shown in Figure~\ref{fig:stormresponseillustration}. Now, it is clear that the truck uses the road network to reach the poles and the power lines of the power grid to do the necessary repairs. In Section~\ref{sec:PG_Network}, we explain in detail the different components of the power grid that are carried by the poles.

Let $T_{tij}$ be the travel time from node $i$ to node $j$ {\color{black} at time $t$ (that is, it is known at time $t$, and random for any decisions made before time $t$)}.  We assume that every point on the power network is accessible by the road network, but the configuration of the road network does not match the power network, as illustrated in Figure 1.

\subsection{The Power Grid}\label{sec:PG_Network}
The power grid consists of substations, poles that carry the protective devices, power lines, and transformers through which customers are fed with power as shown in Figure~\ref{fig:stormresponseillustration}. {\color{black} We ignore buried cables since they are not vulnerable to storm damage. Appendix~\ref{sec:DPG} defines some of the vocabulary used to describe the power grid and additional notation describing the topology of the grid.}

We assume that a storm has blown across the grid generating some faults which result in outages to some customers.   The information that we rely on to restore the grid includes the grid structure,  the phone calls we have received, and the path of the storm.  It is assumed that a fault across a power line includes the fault across its connection ends which could be a transformer or a protective device. This can be justified by the fact that a fault across a power line or one of the components at its end connections results in power outage to the same set of customers in the power system.


Let $\mathcal{U}$ be the set of circuits in the power system, $\mathcal{U}=\{u: u=1,\ldots,U\}$ where each circuit can be represented by a tree that is rooted at the substation. {\color{black} Recalling that $\mathcal{I}=\{i: i=1,\ldots,I\}$ is the set of poles that make up the power grid, let $\mathcal{I}_u=\{i: i=1,\ldots,I_u\}$ be the set of nodes of circuit~$u$ which are mounted on the poles (a node can be either a transformer or a protective device).}  Also, let $\mathfrak{L}_u$ be the set of power lines on circuit~$u$ where power line $i\in \mathfrak{L}_u$ feeds node $i\in \mathcal{I}_u$ with power.


In this paper, we assume the EUCs rely on the calls of the customers that lost power. Let $H_t=\{H_{tui}: \forall i\in \mathcal{I}_u, \forall u\in \mathcal{U}\}$ be a random vector where $H_{tui}$ is a random variable representing the number of phone calls received from node $i$ on circuit~$u$ {\color{black} by time $t$.  We also let $\hat{H}_t$ be the vector of new phone calls that arrived between $t-1$ and $t$, which allows us to write
\begin{eqnarray*}
  H_{t+1} &=& H_t + \hat{H}_{t+1}.
\end{eqnarray*}
}

{\color{black} We next need to determine the power lines that have lost power from the phone calls received by time $t$.} For each circuit~$u$, let $L_{tu}=\{L_{tui}: \forall i\in \mathfrak{L}_u\}$ be a random vector representing the possible realizations of power lines that faulted by time $t$  on circuit~$u$ where $L_{tui}$ is a random variable indicating whether power line $i$ on circuit~$u$ has faulted; we assume $L_{tui}=1$ if the power line faulted and $L_{tui}=0$ otherwise.

{\color{black} Let $R_{ui}$ be a random variable representing the required repair time for power line $i$ on circuit~$u$ which depends on the fault type (which means repair times are random).} Now, we can define $\tau_{tij}=T_{tij} + \sum_{u\in\mathcal{U}}R_{uj}$ as the total time required to go from node $i$ to $j$ which accounts for the travel and repair times of all power lines that feed pole $j$ with power.

Let $\omega$ be a sample realization of the random variables. At time $t$ and according to sample path~$\omega$, let $H_t(\omega)$ be an outcome of customer calls, $L_{tu}(\hspace{-0.03cm}\omega \hspace{-0.05cm})$ be an outcome of the power line faults on circuit~$u$, $R_u(\omega)$ be an outcome of the repair times on circuit~$u$ and $T_{t}(\omega)$ be an outcome of the travel conditions. Then, we can define at time $t$ and according to sample path~$\omega$, $p(H_t=H_t(\omega))$ as the probability of the set of received calls, $p(L_{tu}=L_{tu}(\hspace{-0.03cm}\omega \hspace{-0.05cm}))$ as the prior probability of power line faults on circuit~$u$, $p(R_u=R_u(\omega))$ as the probability of the repair time on circuit~$u$ and  $p(T_t=T_t(\omega))$ as the probability of the travel time.


Define $\Omega$ as the set of all outcomes, $\cal F$ as the set of events and $\cal P$ as the probability measure on $(\Omega,\cal F)$, so that $(\Omega,\cal F,\cal P)$ is the probability space. We have $\Omega\subseteq \mathcal{F} $ formed by a set of scenarios where each scenario~$\omega$ indicates a specific set of calls, a set of power line faults, fault types and travel conditions (traffic). Let $p(\omega)$ be the probability of scenario~$\omega$, where $\sum_{\omega\in\Omega}p(\omega) = 1$.  

{\color{black}
\subsection{Belief State: Probability Model for Power Line Faults}\label{sec:belief_state}
This section discusses the model used to develop a belief about power line faults in the grid. First, the EUC should estimate the prior probability, $p(L_{tui})$, of a fault in power line $i$ on circuit~$u$ at time $t$ based on the storm pattern, the structure of the grid, and the EUC knowledge of the power line conditions and environment as explained in~(\citet{ABP14}). In the worst case, the EUC can assume a uniform probability of fault for all power lines, and then the set of received calls will identify the ones that have most likely faulted.

{\color{black} Let $p(H_t)$ be the probability of the set of lights-out calls $H_t$,} given the probability $\rho$ that a customer calls in to report that their lights are out. \citet{ABP14} derives $p(H_t)$ as a function of $\rho$, the structure of the grid, and the distribution of customers across the grid, and the prior probability of line faults; for this reason, this section simply sketches this logic.
}

For any realization of $L_{tu}$ and $H_t$, the posterior probability of fault of the power lines on circuit~$u$ given the phone calls is calculated using Bayes' theorem as follows:
\begin{eqnarray}
p(L_{tu}|H_t) = \frac{p(H_t|L_{tu})p(L_{tu})}{p(H_t)},
\end{eqnarray}
where $p(H_t|L_{tu})$ is the likelihood of the calls given the power line faults on circuit~$u$ and the expressions have been derived in~(\citet{ABP14}[equation 5]); it is a function of the locations of the calling customers, the call-in probability which refers to the percentage of customers calling to report an outage and the grid structure.  The posterior probability of fault of power line $i$ on circuit~$u$, given the phone calls can be expressed as:
\begin{eqnarray}
p(L_{tui}=1|H_t)=\frac{\sum_{L_{tu}\in \{\mathcal{L}_u\}_{L_{tui}=1}}p(H_t,L_{tu})}{p(H_t)}=\frac{\sum_{L_{tu}\in \{\mathcal{L}_u\}_{L_{tui}=1}}p(H_t|L_{tu})p(L_{tu})}{\sum_{L_{tu}\in\mathcal{L}_u} p(H_t|L_{tu})p(L_{tu})}, \label{eq:eq2}
\end{eqnarray}
\noindent where $\mathcal{L}_u$ is the set containing all power line fault combinations on circuit~$u$ and  $\{\mathcal{L}_u\}_{L_{tui}=1}$ is the set containing a subset of vectors of $\mathcal{L}_u$ where the variable corresponding to power line $i$, i.e., $L_{tui}$, is equal to $1$. Thus, $\{\mathcal{L}_u\}_{L_{tui}=1}$ is the set containing all the combinations of power lines that can fault with power line $i$ on circuit~$u$.

In this work, another factor plays a major role in identifying the fault probability of a power line by time  $t$ which is the trajectory of the truck that is going across the power grid to fix faults. For example, if power line $i$ on circuit~$u$ has been fixed by time $t-1$, then its probability of fault, at time $t$, is $0$ (we assume that once fixed, it does not fault again). 
Let $x_{tij}$ be a binary variable representing whether the utility truck travels from node $i$ to node $j$ using roadway at time $t$. It is assumed that if a truck travels from node~$i$ to node~$j$  at time $t$ then it repairs all the power lines that are attached to pole $j$ if there are faults across them.


Let $x_t$ be a matrix capturing the vector of decisions $(x_{tij})_{i,j\in\mathcal{V}}$ where $x_{tij} = 1$ if the truck is dispatched to $i$ from $j$ at time $t$. Also, let $X_t$ be the trajectory of the truck up to time $t$, i.e., $X_t = (x_{t'})_{t'=0}^t$. The information we are looking for, at time $t$, is the posterior probability, $p(L_{tui}|H_t,X_{t-1})$,  of power line $i$ being in fault given the phone calls and the trajectory of the truck up to time $t-1$ which is given by:
\begin{eqnarray}
&&\hspace{-0.7cm}p(L_{tui}=1|H_t,X_{t-1})=\frac{\sum_{L_{tu}\in \{\mathcal{L}_u\}_{L_{tui}=1}}p(H_t|L_{tu})p(L_{tu}|X_{t-1})}{\sum_{L_{tu}\in\mathcal{L}_u}p(H_t|L_{tu})p(L_{tu}|X_{t-1})}\label{eq:eq3},
\end{eqnarray}
where $p(L_{tu}|X_{t-1})$ is the prior probability of vector $L_{tu}$ being in fault given the route of the truck; the prior probability of the power lines is updated by setting $p(L_{tui}|X_{t-1})=0$ if $\sum_j x_{t'ji}=1, t'\leq t-1$. However, the likelihood $p(H_t|L_{tu})$ is independent of the route of the truck. Accordingly, at each time $t$, the prior probabilities of power line faults should be computed and the set of received calls should account for any new incoming calls between $t-1$ and $t$. After collecting the updated priors and customer calls, the posterior probability of faults is computed using the model described in~\citet{ABP14} upon which the route of the truck is determined.

\section{The Stochastic Optimization Model}~\label{sec:opt_prob}
We now turn our attention to the problem of optimizing the routing of a single utility truck over the transportation network to minimize the amount of time that customers are without power.  Routing decisions have to be made that reflects the geometry of the transportation network, and the beliefs about faults which reflect the prior, lights-out calls from customers, and observations made by the utility trucks.  As the truck moves down street segments where a line has been damaged, the truck can repair the damage (but this takes time).   This is the first rigorous formulation of an information-collecting vehicle routing problem, and we propose a stochastic lookahead policy using Monte Carlo tree search, which appears to be its first formal use in a stochastic vehicle routing problem.

The problem introduces a number of challenges. First, the number of customers whose power is restored after a truck visits a location (even if a repair occurs) is a random variable since the actual value depends also on upstream and downstream outages which are uncertain. Second, each time a truck travels a segment, the information collected (e.g. that there is an outage on that segment, or not) is used to update the probability of outages on all lines. Third, new phone calls are arriving over time, which allows us to update probabilities of outages.  At the same time, as trucks identify and fix outages, other phone calls may become irrelevant. Finally, the time required for a truck to traverse a segment depends to a large extent on whether it finds an outage, and the time required to repair the outage.

\begin{figure}[t!]
\vspace{-0cm}
\centering
\includegraphics[width=2.5in]{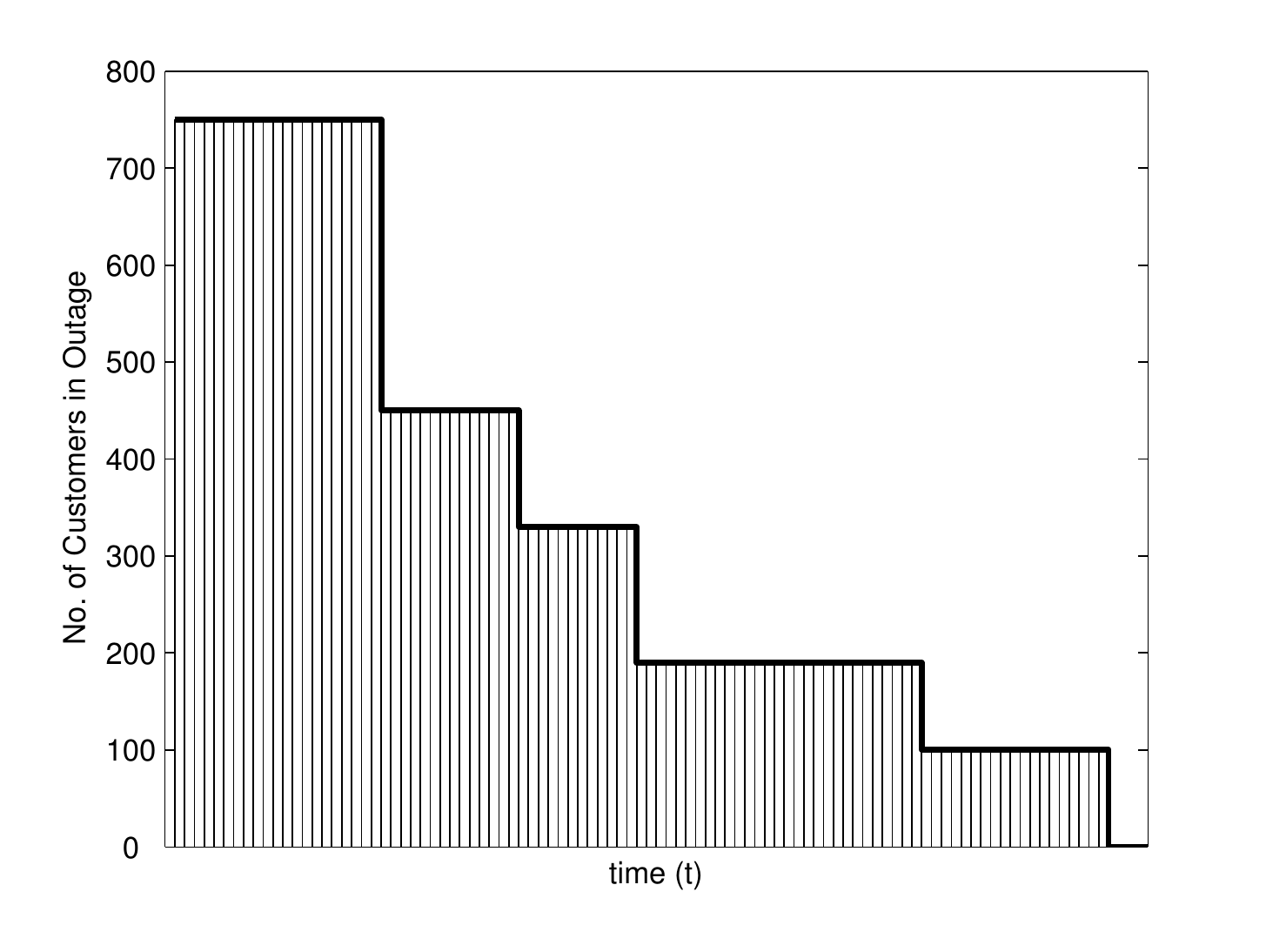}
\vspace{-0.3cm}
\caption{\label{fig:obj}Objective function; outage-minute is represented by the shaded area under the curve.}\vspace{-0.3cm}
\end{figure}

We present the stochastic optimization problem using the canonical modeling framework that consists of states, decisions, exogenous information, transition function and objective function. {\color{black} We index all variables by time $t$, but time will only be measured at instants where the vehicle arrives at a node and has to make a decision.  To simplify notation, time $t$ will represent a counter, so $t+1$ represents the next decision instant.  Also, any variable indexed by $t$ is assumed to be known at time $t$.}
\begin{itemize}
\item State $S_t$ - The information capturing what we know at time $t$. {\color{black} In this work, $S_t=\left(R_t,I_t,B_t\right)$ where $R_t$ represents the physical state that indicates the location of the truck at time $t$, $I_t = H_t$ is other deterministic information which captures the set of received calls by time~$t$, and $B_t$ is the belief state about the faults in the grid discussed in section~\ref{sec:belief_state}. The belief state is given by $B_t=(P_t^L)$ where  $P_t^L$ is a vector containing the probability of all power line faults (i.e., its entries are $p(L_{tui}=1|H_{t})$).} Thus, given the prior probabilities of fault and the set of calls, we can calculate the posterior probabilities of faults which also represent the state of the network at time~$t$.
\item Decision $x_t$ - The vector $x_t = (x_{tij})_{i,j}$  describes the movement of the truck from one location on the grid (represented by ``poles'' that hold transformers or circuit breakers) to the next.  The decision $x_{tij} = 1$ means the truck is moving from pole $i$ to pole $j$ for $i,j \in {\cal I}$.   Let $\mathcal{X}_t$ be the set of poles that a truck in state $S_t$ (which includes the current location) can move to at time $t$. Ideally, $\mathcal{X}_t$ includes movements to all locations in the grid but since the truck is moving on a road network, we can limit $\mathcal{X}_t$ to movements to locations that are within a limited distance from the location of the truck at time $t$. Let $X^\pi(S_t)$ be the policy that determines $x_t \in \mathcal{X}_t$ given~$S_t$.
\item Exogenous information $W_t$ - The new information that arrives between $t-1$ and $t$.  This includes new phone calls ${\hat H}_t$, as well as information about outages (discovered by the utility trucks) and travel (or repair) times.  We denote the exogenous information process by $W_1, W_2, ..., W_T$ where $W_t$ depends on both the state $S_t$ (since it depends on where the truck is located) and the decision~$x_t$.  Let $\omega$ be a sample realization of the information.  For this reason, we let $\Omega^\pi$ be the set of outcomes which depends on the policy $\pi$ that we use to dispatch trucks. 

\item {\color{black} The transition function $S_{t+1}=S^M(S_t,x_t,W_{t+1})=(R_{t+1},I_{t+1},B_{t+1})$ which represents the evolution of the physical, informational, and belief states. The transition function includes all the equations that govern the dynamics of the system, whether it is moving the truck from one location to the next, discovering that a tree has fallen down on a power line, observing new lights-out calls, or using Bayes' theorem to update beliefs.}
\item {\color{black} The objective function - Our objective is to find the policy that maximizes the number of customers with restored power over time which is equivalent to minimizing the number of customers in outage. We define the cost function $C(S_t,x_t)$ as
    \begin{eqnarray*}
    C(S_t,x_t) &=& \parbox[t]{3.5in}{the number of customers in outage at time $t$.}
    \end{eqnarray*}
    $C(S_t,x_t)$ is shown in Figure~\ref{fig:obj} and the sum over time, $\sum_{t=0}^TC(S_t,x_t)$, evaluates the shaded area under the curve which we refer to as ``customer outage-minutes.''
    The total objective while using the policy and following sample path $\omega\in\Omega^\pi$ can be represented as
    \begin{eqnarray*}
    F^\pi(\omega) = \sum_{t=0}^TC\big((S_t(\omega),X^\pi(S_t(\omega))\big).
    \end{eqnarray*}
    Our challenge is to find the policy that solves:
    \begin{eqnarray}
    \min_\pi \mathbb{E}^\pi\left[\sum_{t=0}^TC(S_t,X^{\pi}(S_t))|S_0\right]\label{eq:optimal_policy}
    \end{eqnarray}
    where the expectation is over all possible sequences $W_1,W_2,\ldots,W_T$, which  depend on the decisions taken, and over the belief in $S_0$, which captures all deterministic parameters as well as the priors on where outages might be (based, for example, on storm information and past history).}
\end{itemize}

We refer to this model as the {\it base model} since we are designing policies to solve this problem, captured in the objective function \eqref{eq:optimal_policy}.  Later we are going to introduce the idea of a {\it lookahead model} which is an approximation of the base model.

The objective function captures the amount of time that customers have lost power.  The number of customers whose power is restored due to decision $x_t$ (which moves the truck along a link of the network) is a random variable given the uncertainty about the location of outages, and the location of circuit breakers that open when a line failure occurs. Fixing an outage at one link can allow an upstream breaker to be closed, which turns on the lights of all downstream customers.

The actual ``customer outage-minutes'' objective is represented in Figure~\ref{fig:obj}, but since we cannot {\color{black}determine} the exact number of customers with restored power at each state $S_t$, we evaluate the expected customer outage-minutes. At state $S_t$, the fault probabilities are updated according to (\ref{eq:eq3}) and the expected number of customers in outage~is:
    \begin{eqnarray*}
    C(S_t,x_t)=   \sum_{u\in\mathcal{U}}\sum_{s\in S_u} \left(1-\prod_{k\in Q_s}p(L_{tuk}=0)\right)\sum_{k\in s} n_{uk},
 \label{eq:ecom}
    \end{eqnarray*}
where $Q_s$ represents the set of power lines that if a fault results in an outage to segment $s$, the term in parenthesis is the probability of at least one fault across $Q_s$ that causes $s$ to be in outage and $\sum_{k\in s} n_{uk}$ is the number of customers across segment $s$.

So far, we have defined the fundamental elements influencing truck routing. In this work, the flow of information and state transition occurs in the following sequence:
\begin{eqnarray}
(S_0,x_0,S^x_0,W_1,S_1,x_1,S^x_1,W_2,S_2,\ldots,S_t,x_t,S^x_t,W_t,\ldots,S_T)
\end{eqnarray}
where $S^x_t$ represents the state after taking decision $x_t$, known as the post-decision state. Based on $x_t$, we see a specific exogenous information $W_{t+1}$ resulting from the decision $x_t$.

The distinguishing feature of this problem relative to the classical stochastic vehicle routing literature is the dimension that the utility truck is collecting information, and that our state variable includes the physical state of the truck, our current state of knowledge concerning the probability of outages, and the history of phone calls (there is a reason we have to retain the history).  A decision to dispatch a truck from $i$ to $j$ has to consider not only the change in the physical state of the truck but also the value of the information that is collected while traversing from $i$ to $j$.

\section{Designing policies}
\label{sec:designingpolicies}
{\color{black}
\citet{Powell2022RLSO} describes two strategies for solving the optimization in equation \eqref{eq:optimal_policy}, each of which can be divided into two classes.  These are:
\begin{itemize}
\item The policy search class - These are policies that do well when simulated in our objective function \eqref{eq:optimal_policy}.  We may search over classes of policies, or we may search within a class by tuning parameters.  We can divide this class into two subclasses:
    \begin{itemize}
    \item Policy function approximations (PFAs) - These are analytic functions that map states to actions.  They may be lookup tables (when at this node turn left), linear models (often called ``affine policies'') such as
        \begin{eqnarray*}
        X^{PFA}(S_t|\theta) &=& \sum_{f\in{\cal F}} \theta_f \phi_f(S_t).
        \end{eqnarray*}
        or nonlinear models such as a logistic regression or neural network.
    \item Cost function approximations (CFAs) - These are parameterized optimization problems.  For example, we might optimize the assignment of resources to tasks with bonuses and penalties to achieve different objectives.  We can write these generally as
        \begin{eqnarray*}
        X^{CFA}(S_t|\theta) = \argmin_{x_t\in{\cal X}^\pi(\theta)} {\bar C}^\pi(S_t,x_t|\theta).
        \end{eqnarray*}
    \end{itemize}
    Note that both PFAs and CFAs require tunable parameters $\theta$ that would have to be tuned using the objective function in \eqref{eq:optimal_policy}.
\item The lookahead class - These are policies that depend on approximations of the downstream impact of making a decision now.  These can also be divided into two classes:
    \begin{itemize}
    \item Policies based on value function approximations (VFAs) - These are the policies based on Bellman's equation such as
        \begin{eqnarray*}
        X^{VFA}(S_t) = \argmin_{x_t\in{\cal X}_t} \big(C(S_t,x_t) + {\bar V}^x_t(S_t,x_t)\big),
        \end{eqnarray*}
        where ${\bar V}^x_t$ is the post-decision value function approximation that depends on the post-decision state $S^x_t$ which depends (deterministically) on $S_t$ and $x_t$.
    \item Direct lookahead policies (DLAs) - These are policies computed by optimizing over some horizon.  The most familiar are deterministic lookaheads (think of navigation systems for cars).  We describe the full DLA below.
    \end{itemize}
\end{itemize}

The utility currently uses a rule-based policy, which falls in the PFA class (``send a truck to the next caller reporting an outage'').  A common method for routing vehicles is to solve a deterministic lookahead, but introduces tunable parameters (such as schedule slack) to handle uncertainties (this would be a hybrid CFA/DLA).

In this paper, we use a policy based on a stochastic lookahead~model that simultaneously plans both the physical movement of the truck and beliefs about outages.  Lookahead policies are almost universally used in vehicle routing since the right decision now depends on what we are going to do with the vehicle in the future.  However, there are choices to be made in how we model the future.

An optimal lookahead model solves the base model by making a decision $x_t$ given that we are in state $S_t$ that optimizes all remaining costs into the future.  This can be written as
}
\begin{eqnarray}
X^*_t(S_t)=\arg\min_{x_t\in{\mathcal{X}_t(S_t)}}\left(C(S_t,x_t)+\mathbb{E} \left\{\min_\pi \mathbb{E}\left\{ \sum_{t'=t+1}^T C(S_{t'},X^\pi_{t'}(S_{t'}))|S_{t+1}\right\}|S_t,x_t\right\}\right), \label{eq:opt_policy}
\end{eqnarray}
where $S_{t+1}=S^M\left(S_t, x_t,W_{t+1}\right)$.

In practice, computing (\ref{eq:opt_policy}) is compuftationally intractable, requiring that we introduce an approximation. There are six strategies that are typically used to simplify lookahead models:
\begin{itemize}
\item[1) ] Limiting the horizon by reducing it from $(t,T)$ to $(t,t+H)$, where we are repurposing ``$H$'' to planning horizon.
\item[2) ] Discretizing the time, states and decisions to make the model computationally tractable.
\item[3) ] Aggregating the outcome or sampling by using Monte Carlo sampling to choose a small set, $\tilde{\Omega}_t$, of possible outcomes between $t$ and $t+H$.
\item[4)] Stage aggregation which represents the process of revealing information before making another decision (see \citet{BL11} for a thorough introduction to stochastic programming using scenario trees). A common approximation is a two-stage formulation, where we make
a decision $x_t$, then observe all future events (until $t+H$), and finally make all remaining decisions. 
\item[5)] {\color{black} Policy simplification - Stochastic lookahead policies are, themselves, stochastic optimization problems that require a ``policy within the policy.'' It is not unusual to substitute a simplified policy (sometimes called a rollout policy, but this is not the only choice) for the ``policy within the policy.''}
\item[6)] Dimensionality reduction where we ignore some variables in our lookahead model as a form of simplification. For example, a forecast of future incoming phone calls can add a number of dimensions to the state variable. While we have to track these in the original model, we can hold them fixed in the lookahead model, and then ignore them in the state variable (these become latent variables).
\end{itemize}

{\color{black} It is not unusual for lookahead models to incorporate all six of these forms of simplifications.  In our model, we only made one simplification: we do not model new ``lights-out'' calls in the lookahead model, given by $\hat{H}_t$. Rather, the only source of new information we captured in the lookahead model was whether a tree had fallen down on a link.  The reason for the simplification is that the Bayesian updating for the lights-out calls is relatively expensive, and is magnified inside the lookahead model because we might do this hundreds of times in a single call to the MCTS algorithm. In principle we could use the exact same model in the lookahead as we use in the base model, but while possible, it is computationally expensive.}


All variables in the lookahead model are indexed by $(t,t')$ where $t$ represents when the decision is being made (which fixes the information content) while $t'$ is the time within the lookahead model. We also use tilde's to avoid confusion between the lookahead model (which often uses a variety of approximations) and the real model. Thus states, decisions and exogenous information are represented as $\tilde{S}_{tt'}$, $\tilde{x}_{tt'}$ and $\tilde{W}_{tt'}$, where states are updated with the transition function $\tilde{S}^M(\tilde{S}_{tt'},\tilde{x}_{tt'},\tilde{W}_{t,t'+1})$.  Using this notation,  the stochastic lookahead policy in \eqref{eq:opt_policy} becomes
\begin{eqnarray}
X^{DLA}_t(S_t) &=& \argmin_{x_t} \left(C(S_t,x_t) + {\tilde E} \left\{\min_{{\tilde \pi}} {\tilde E} \left\{\sum_{t'=t+1}^T C({\tilde S}_{tt'},{\tilde X}^{\tilde \pi}({\tilde S}_{tt'})) | {\tilde S}_{t,t+1}\right\}|S_t,x_t\right\}\right). \label{eq:optimalpolicyLAapproxpi}
\end{eqnarray}

{\color{black}
Writing the approximate lookahead policy in the form given in equation \eqref{eq:optimalpolicyLAapproxpi} creates a natural bridge to policies that depend on the idea known as a rollout policy.
\cite{Ulmer2019} proposes an ``offline-online'' method that combines value function approximations (which are estimated offline) with a rollout policy that is performed dynamically as a vehicle moves through the system (this is the online component).  The offline VFA produces a value function approximation ${\bar V}^x_t(S^x_t)$ around the post-decision state $S^x_t$ (which is a deterministic function of $S_t$ and $x_t$) that produces a rollout policy
\begin{eqnarray}
X^{rollout}_t(S_t) = \argmin_x \big(C(S_t,x) + {\bar V}^x_t(S^x_t)\big). \label{eq:rolloutvfa}
\end{eqnarray}
We first note that the rollout policy in \eqref{eq:rolloutvfa} represents the solution of ``$\min_{\tilde \pi} {\tilde E}$'' in \eqref{eq:optimalpolicyLAapproxpi} by using a VFA-based policy. This rollout policy is then simulated over some horizon to obtain a sampled estimate of the value of being in an initial state $S_t$.  This is likely to be more accurate than just building an approximate VFA ${\bar V}_t(S_t)$ of the value of being in state $S_t$ since value function approximations are limited by what we can represent using a statistical model, whereas a rollout policy can capture much more complex state dynamics.  As is pointed out in \cite{Ulmer2019}, the VFA-based rollout policy captures the overall structure (and avoids the need to construct ad-hoc policies), while the rollout policy captures local details of the problem.



MCTS takes the challenge of designing a lookahead policy ${\tilde \pi}$ a step further by actually modeling the full nesting of decisions in the future (that is, modeling the full decision tree). To do this we can rewrite the optimal policy in equation \eqref{eq:opt_policy} by explicit writing the sequencing of decision $x_t$, information $W_{t+1}$, followed by the decision $x_{t+1}$, and so on. Thus, \eqref{eq:opt_policy} is equivalent to:
\begin{eqnarray}
X^*_t(S_t)&=&\arg\min_{x_t\in{\mathcal{X}_t(S_t)}}\bigg(C(S_t,x_t)+\mathbb{E}_{W_{t+1}}\bigg[\min_{x_{t+1}\in{\mathcal{X}_{t+1}}}C(S_{t+1},x_{t+1})+\mathbb{E}_{W_{t+2}}\bigg[ \ldots + \nonumber\\
&&\hspace{5.5cm}\mathbb{E}_{W_{T}}\bigg[C(S_{T})|S^x_{T-1}\bigg]\ldots\bigg]|S^x_{t+1}\bigg]|S_t,x_t \bigg]\bigg).\label{eq:eqv_opt_policy}
\end{eqnarray}
We then make the transition from the full lookahead model to our approximate lookahead model which produces
\begin{eqnarray}
\hspace{-0.5cm}X^{DLA}_t(S_t)&=&\arg\min_{x_t\in{\mathcal{X}_t(S_t)}}\bigg(C(S_t,x_t)+\tilde{\mathbb{E}}_{\tilde{W}_{t,t+1}}\bigg[\min_{\tilde{x}_{t,t+1}\in{\mathcal{\tilde{X}}_{t,t+1}(\tilde{S}_{t,t+1})}}\tilde{C}(\tilde{S}_{t,t+1},\tilde{x}_{t,t+1})+ \nonumber\\
&&\hspace{-0.7cm}\mathbb{\tilde{E}}_{\tilde{W}_{t,t+2}}\bigg[ \ldots\mathbb{\tilde{E}}_{\tilde{W}_{t,t+H}}\bigg[\tilde{C}(\tilde{S}_{t,t+H})|\tilde{S}^x_{t,t+H-1}\bigg]\ldots\bigg]|\tilde{S}^x_{t,t+1}\bigg]|S_t,x_t \bigg]\bigg),\label{eq:alo}
\end{eqnarray}
where the expectation $\tilde{\mathbb{E}}\{.|S_t,x_t\}$ is over the sample space in $\tilde{\Omega}_{t,t+1}$ which is constructed given that we are in state $S_t$ at time $t$. We emphasize that the state captures the physical state (location of the vehicle, known grid outages) and the belief state, consisting of probabilities of outages for portions we have not visited. When computing this policy, we start in a particular state $S_t$, but
then step forward in time using:
\begin{eqnarray}
\tilde{S}_{t,t'+1}=\tilde{S}^M(\tilde{S}_{tt'},\tilde{x}_{tt'},\tilde{W}_{t,t'+1}), ~~t'=t,\ldots,t+H-1.
\end{eqnarray}
}

Generating the entire tree is computationally intractable.  For this reason, we turn to a popular strategy developed in the computer science community known as Monte Carlo tree search (MCTS)~ (\citet{Munos14,BPW12,CHS08,KS06}) which uses an intelligent sampling procedure to create a partial tree that can be solved, and which asymptotically ensures that we would find the optimal solutio of the lookahead model with a large-enough search budget.  Given that we are solving the one-truck problem, we can formulate the lookahead model as a decision tree, taking advantage of the property that the number of possible decisions for a truck at any point in time is reasonaly small, since it is limited to the number of locations on the distribution grid that a truck can drive to next.  Further, if we limit the random information to whether the truck finds a fault or not, then the random variables are binomial.  However, even with these restrictions, a decision tree will still grow exponentially. 


We propose  MCTS as the look-ahead policy to solve the original problem in the following way. Given a current state, $S_t$, which depends on the location of the truck and the probability of faults at time~$t$,  MCTS should decide where to move the truck next in order to minimize the objective represented by the customer outage-minutes. So, starting from the current state $S_t$ as the root node, MCTS successively builds the look-ahead tree over the state-space. Finally, the move that corresponds to the highest value from the root node will be taken.

At this stage, the belief state is updated taking into account any incoming exogenous information such as whether there was a fault or not by taking the move that was the outcome of the previous step, the newly arrived phone calls of customers, and the consumed travel/repair times. Then, the whole process is repeated until power is restored to the whole power grid based on the values of the probability model. The pseudo-code of the proposed lookahead policy to solve the utility truck routing problem is presented in Algorithm~\ref{alg:adp}.

We note that the belief state is high dimensional and continuous, but MCTS enjoys the property that it is not sensitive to the size of the state space, but its complexity grows quickly with the number of possible decisions that can be made at each point.


\begin{algorithm}[t!]{\footnotesize
\caption{Lookahead Policy for Utility Truck Routing}
\label{alg:adp}
\begin{algorithmic}[t!]
\STATE\textbf{Step 0.} \textbf{Initialization}:  Initialize the state $S_0=(R_0,B_0)$ where $B_0=(P_0^L,H_0)$. Set $t\leftarrow0$.
\STATE\textbf{Step 1.}~\textbf{While} $\left(p(L_{tui}|H_t,X_{t-1})\geq \epsilon^{thr}, \forall i,\forall u\right)$ do:
\STATE \hspace{1.2cm}\textbf{Step 1a.} At time $t$, fix the set of calls $H_t$ and determine $\tilde{\Omega}_t$.
\STATE \hspace{1.2cm}\textbf{Step 1b.} Call $MCTS(S_t)$ (see section \ref{sec:MCTS})  to solve (\ref{eq:alo}) that determines $x^*_t$.
\STATE \hspace{1.2cm}\textbf{Step 1c.} Set $x_t\leftarrow x^*_t$ and move the truck according to $x_t$.
\STATE \hspace{1.2cm}\textbf{Step 1d.} Update $S_{t+1}=S^M(S_t,x_t,W_{t+1})=(R_{t+1},I_{t+1},B_{t+1})$ until the truck reaches the destination node \STATE \hspace{2.5cm}  set by $x^*_t$, say at time $t'$.
\STATE \hspace{1.2cm}\textbf{Step 1e.} $t\leftarrow t'$
\STATE \hspace{1cm}\textbf{End While}
\end{algorithmic}\vspace{-0.0cm}
}
\end{algorithm}

{\color{black}

\section{Optimistic Monte Carlo Tree Search}\label{sec:MCTS}

In this section we provide a detailed summary of optimistic MCTS.  We begin with a general overview of MCTS, followed by a detailed discussion of the four steps that make up MCTS.  We close with a discussion of the convergence of MCTS using the concept of information relaxation, which is the basis of optimistic MCTS.

\subsection{The general strategy}
MCTS builds a search tree until we hit a limit on the number of iterations.  Each iteration of MCTS is formed of four steps: Selection, Expansion, Simulation, and Backpropagation~(\citet{CHS08}). The \emph{Selection} step involves selecting a decision successively starting from the initial state until an expandable state is reached.  The \emph{Expansion} step adds one or more states to the tree. The \emph{Simulation} step is referred to as the ``simulation policy'' (or ``rollout policy'') which is used to evaluate the value of the newly added state. Finally, in the \emph{Backpropagation} step, the value of the newly added state is backpropagated to update the value functions of all predecessor states.

MCTS is a forward algorithm that searches forward in time, with heuristic methods for approximating the value of states that have not yet been reached.  MCTS uses the simulation policy to obtain initial estimates of states that have not been visited before.  This is typically done by using the structure of the problem to design a reasonable heuristic \citep{BPW12}. \citet{Goodson2017} uses rollout policies tuned for vehicle routing problems that work quite well.

Simulation policies are heuristics that overestimate (on average) the cost of being in a state.  We refer to this as ``pessimistic MCTS.''  In our paper, we use a novel approach that samples the future, and then solves a deterministic optimization problem that underestimates the cost of being in a state since decisions are allowed to see events into the future.  This strategy can be viewed as a form of information relaxation (see \cite{Brown2010}).  We refer to this as ``optimistic MCTS''; it is computationally more expensive, but produces more rapid convergence (in terms of iterations).  Perhaps more important is that it avoids the need to design ad-hoc rollout policies.  It can also serve as a useful benchmark for evaluating heuristic rollout policies.

It has been shown (see \cite{Kocsis2006e} and \citet{ACT13}) that pessimistic MCTS is asymptotically optimal. Building on our work, \cite{JiangMCTS2020} shows that optimistic MCTS is also asymptotically optimal.  However, the mechanisms of the two proofs are quite different.  Pessimistic MCTS depends on the use of upper confidence bounding to force the algorithm to explore each decision infinitely often (in the limit). Optimistic MCTS does not require this, depending instead on the ability to prune actions that do not appear attractive basic on optimistic estimates.


The differences in these proofs has practical ramifications.  Pessimistic MCTS is sensitive to the number of actions.  For example, it would handle well an action space for a truck that consists only of whether the truck should go straight, right or left at an intersection.  In our problem, however, trucks choose which point on the grid to visit next, after which we solve a shortest path problem to determine the path over the transportation network.  As a result, a truck has a large number of choices.  Classical pessimistic MCTS would require sampling all of these.


One feature of MCTS that is critical for our application is that it is relatively insensitive to the size and complexity of the state space.  This is important because the belief state is both continuous and very high dimensional.  Recall that the belief state $B_t = (P^L_{tij})_{i,j\in\mathcal{V}}$ which means it is a vector dimensioned by the number of segments in our grid.  We do have to recognize that we have to store the state variable for each node in the MCTS tree, which could become problematic if we wanted to capture conditional probabilities.
}



In stochastic MCTS, which is also referred to as sampled MCTS, the states in the tree are associated with the following~data:
 \begin{itemize}
 \item[1) ] The pre-decision value function, $\tilde{V}_{tt'}(\tilde{S}_{tt'})$, the post-decision value function, $\tilde{V}^x_{tt'}(\tilde{S}^x_{tt'})$, and cost, $\tilde{C}(\tilde{S}_{tt'},\tilde{x}_{tt'})$, which represent the value function and cost of being in state $\tilde{S}_{tt'}$ and taking decision $\tilde{x}_{tt'}$, respectively.
 \item[2) ] The visit count, $N(\tilde{S}_{tt'})$, which represents the number of rollouts that included state $\tilde{S}_{tt'}$.
 \item[3) ] The count of the state-decision, $N(\tilde{S}_{tt'},\tilde{x}_{tt'})$, which represents the number of times decision $\tilde{x}_{tt'}$ was taken from state $\tilde{S}_{tt'}$.
 \item[4) ] The set of decisions, $\tilde{\mathcal{X}}_{tt'}(\tilde{S}_{tt'})$, and set of possible random outcomes, $\tilde{\Omega}_{t,t'+1}(\tilde{S}_{tt'}^{x})$; $\tilde{\mathcal{X}}_{tt'}(\tilde{S}_{tt'})$ is the set of decisions that the truck would face moving over a road network given that it is at state $\tilde{S}_{tt'}$. For state $\tilde{S}_{tt'}$, let $\tilde{\mathcal{X}}_{tt'}^e(\tilde{S}_{tt'})$ be the set of decisions that has been explored by the truck (i.e., expanded in the tree)  by time $t'$ and let $\tilde{\mathcal{X}}_{tt'}^u(\tilde{S}_{tt'})$ be its complement set which represents the set of unexplored decisions in the tree by time $t'$. Similarly,  $\tilde{\Omega}_{t,t'+1}(\tilde{S}_{tt'}^{x})$ is the set of all possible random events that can take place at time $t'+1$ given state $\tilde{S}_{tt'}^{x}$, $\tilde{\Omega}_{t,t'+1}^e(\tilde{S}_{tt'}^{x})$ is the set of explored events and $\tilde{\Omega}_{t,t'+1}^u(\tilde{S}_{tt'}^{x})$ is its complement.
 \end{itemize}

\begin{figure}[t!]
\vspace{-0cm}
\centering
\includegraphics[scale=0.36]{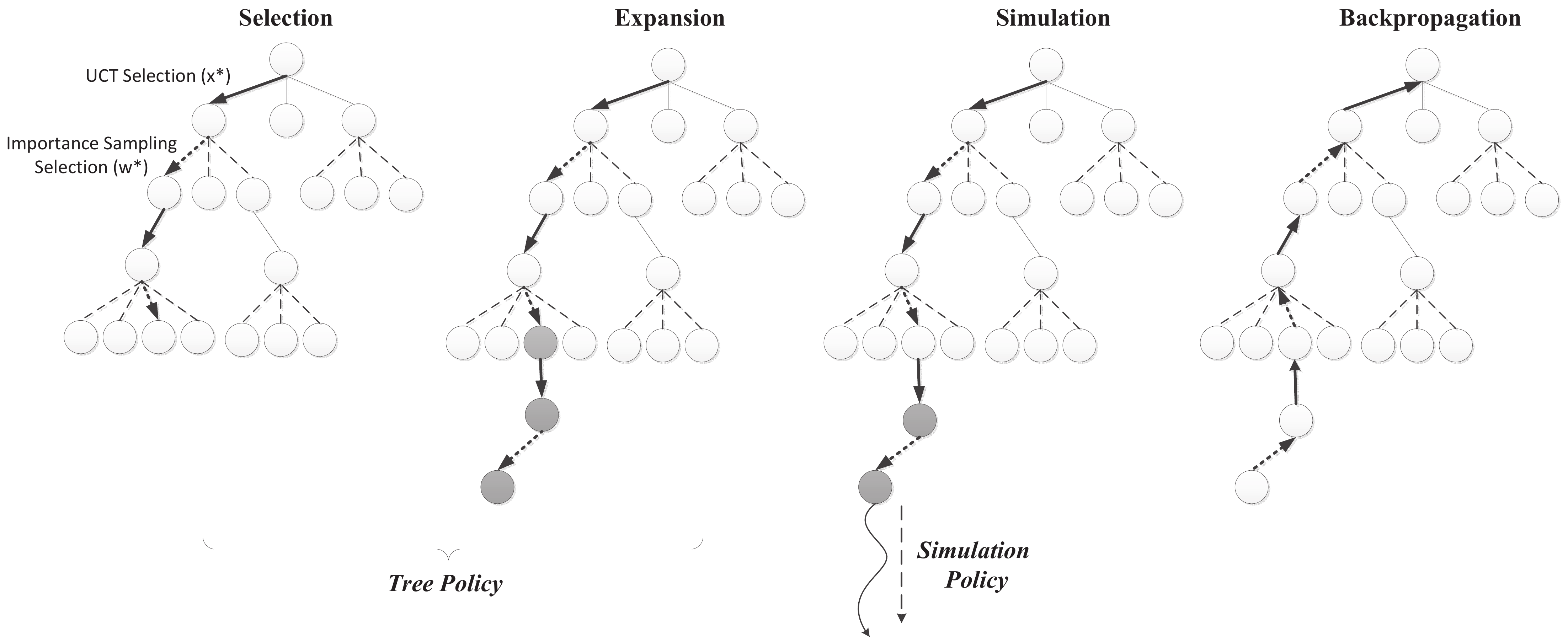}
\vspace{-0.3cm}
\caption{\label{fig:S_MCTS}One iteration of the proposed optimistic MCTS.}\vspace{-0.3cm}
\end{figure}

In sampled MCTS, assume that we are at a \emph{pre-decision state} $\tilde{S}_{tt'}$ and decide to take decision $\tilde{x}_{tt'}$ which takes us to the \emph{post-decision state} $\tilde{S}_{tt'}^{x}$. In real life, while the truck is moving to the location specified by $x_{t}$, it encounters first the travel time which depends on traffic. Then,  it discovers the fault type at the intended location which determines the repair time and affects the number of customers with restored power. In the lookahead model, upon determining $\tilde{x}_{tt'}$, we sample one of the possible exogenous realizations $\tilde{W}_{t,t'+1}$ which immediately informs us about the expected travel and repair times and thus we can immediately know the time, $\tau(\tilde{x}_{tt'},\tilde{W}_{t,t'+1})$, required by the truck to arrive to the destination node specified by $\tilde{x}_{tt'}$. Thus, in MCTS, we define the stochastic transition function as $\tilde{S}_{t,t'+\tau(\tilde{x}_{tt'},\tilde{W}_{t,t'+1})}=\tilde{S}^{M,x}(\tilde{S}_{tt'}^{x},\tilde{W}_{t,t'+1})$. It is obvious now that taking the same decision $\tilde{x}_{tt'}$ from the same state $\tilde{S}_{tt'}$ results in a different outcome state based on the exogenous information $\tilde{W}_{t,t'+1}$.


We assume that MCTS has a computational budget of $n^{iter}$ iterations. The steps of MCTS can be grouped into two main policies: a \emph{Tree Policy } (formed of the \emph{Selection} and \emph{Expansion} steps) and a \emph{Simulation Policy} (formed of the \emph{Simulation} step) as shown in Figure~\ref{fig:S_MCTS}. After terminating the tree search, then the decision that corresponds to best value from the root node is chosen. The pseudo-code for the optimistic MCTS is presented in Algorithm~\ref{alg:MCTS}.




\begin{algorithm}[t!] \footnotesize
\caption{Optimistic MCTS}
\label{alg:MCTS}
\begin{algorithmic}[t!]
\STATE \textbf{function} $MCTS(S_t)$
\STATE \hspace{0.5cm}Create root node $\tilde{S}_{tt}$ with state $S_t$; set iteration counter $n=0$
\STATE \hspace{0.5cm}{\bf while} $n< n^{iter}$
\STATE \hspace{1cm} $\tilde{S}_{tt'}\leftarrow TreePolicy(\tilde{S}_{tt})$
\STATE \hspace{1cm} $\tilde{V}_{tt'}(\tilde{S}_{tt'})\leftarrow SimPolicy(\tilde{S}_{tt'})$
\STATE \hspace{1cm} $Backup(\tilde{S}_{tt'},\tilde{V}_{tt'}(\tilde{S}_{tt'}))$
\STATE \hspace{1cm} $n\leftarrow n+1$
\STATE \hspace{0.5cm}{\bf end while}
\STATE \textbf{return} $x_{t}^*=\arg\min_{\tilde{x}_{tt}\in \tilde{\mathcal{X}}_{tt}^e(\tilde{S}_{tt})}\tilde{C}(\tilde{S}_{tt},\tilde{x}_{tt})+\tilde{V}^x_{tt}(\tilde{S}^x_{tt})$
\STATE
\STATE \textbf{function} $TreePolicy(\tilde{S}_{tt})$
\STATE \hspace{0.5cm} $t'\leftarrow t$
\STATE \textbf{while} $\tilde{S}_{tt'}$  is non-terminal \textbf{do}
\STATE \hspace{0.2cm} \textbf{if} $|\tilde{\mathcal{X}}_{tt'}^e(\tilde{S}_{tt'})|<d^{thr}$ \textbf{do}(Expanding a decision out of a pre-decision state)
\STATE \hspace{0.5cm} choose decision  ${\tilde{x}_{tt'}}^*$ optimistically by using a two-stage lookahead model
\STATE \hspace{0.5cm} $\tilde{S}_{tt'}^{x}=\tilde{S}^M(\tilde{S}_{tt'},{\tilde{x}_{tt'}}^*)$ (Expansion step)
\STATE \hspace{0.5cm}  $\tilde{\mathcal{X}}_{tt'}^e(\tilde{S}_{tt'})\leftarrow \tilde{\mathcal{X}}_{tt'}^e(\tilde{S}_{tt'})\bigcup\{\tilde{x}_{tt'}^*\}$
\STATE \hspace{0.5cm}  $\tilde{\mathcal{X}}_{tt'}^u(\tilde{S}_{tt'})\leftarrow \tilde{\mathcal{X}}_{tt'}^u(\tilde{S}_{tt'})-\{\tilde{x}_{tt'}^*\}$
\STATE \hspace{0.2cm} \textbf{else}
\STATE \hspace{0.5cm} $\tilde{x}_{tt'}^*=\arg\min_{\tilde{x}_{tt'}\in \tilde{\mathcal{X}}_{tt'}^e(\tilde{S}_{tt'})}\left(-\left(\tilde{C}(\tilde{S}_{tt'},\tilde{x}_{tt'})+\tilde{V}^x_{tt'}(\tilde{S}^x_{tt'})\right)+\alpha\sqrt{\frac{\ln N(\tilde{S}_{tt'})}{N(\tilde{S}_{tt'},\tilde{x}_{tt'})}}\right)$
\STATE \hspace{0.5cm} $\tilde{S}_{tt'}^{x}=\tilde{S}^M(\tilde{S}_{tt'},{\tilde{x}_{tt'}}^*)$
\STATE \hspace{0.2cm} \textbf{end if}
\STATE \hspace{0.2cm} \textbf{if} $|\tilde{\Omega}_{t,t'+1}^e(\tilde{S}_{tt'}^{x})|<e^{thr}$  \textbf{do} (Expanding an exogenous outcome out of a post-decision state)
\STATE \hspace{0.5cm} choose exogenous event ${\tilde{W}_{t,t'+1}}$ according to importance sampling with uniform distribution $g(\tilde{W}_{t,t'+1})$,\STATE \hspace{13cm} $\forall \tilde{\omega}\in\tilde{\Omega}_{t,t'+1}^u(\tilde{S}_{tt'}^{x})$
\STATE \hspace{0.5cm} $\tilde{S}_{t,t'+\tau(\tilde{x}_{tt'},\tilde{W}_{t,t'+1})}=\tilde{S}^{M,x}(\tilde{S}_{tt'}^{x},\tilde{W}_{t,t'+1})$ (Expansion step)
\STATE \hspace{0.5cm} $\tilde{\Omega}_{t,t'+1}^e(\tilde{S}_{tt'}^{x}) \leftarrow\tilde{\Omega}_{t,t'+1}^e(\tilde{S}_{tt'}^{x})\bigcup\{\tilde{W}_{t,t'+1}\}$
\STATE \hspace{0.5cm} $\tilde{\Omega}_{t,t'+1}^u(\tilde{S}_{tt'}^{x}) \leftarrow\tilde{\Omega}_{t,t'+1}^u(\tilde{S}_{tt'}^{x})-\{\tilde{W}_{t,t'+1}\}$
\STATE \hspace{0.5cm} $t'\leftarrow t'+\tau(\tilde{x}_{tt'},\tilde{W}_{t,t'+1})$
\STATE \hspace{0.5cm} \textbf{return} $\tilde{S}_{tt'}$ (stops execution of \textbf{while} loop)
\STATE \hspace{0.2cm} \textbf{else}
\STATE \hspace{0.5cm} choose exogenous event ${\tilde{W}_{t,t'+1}}$ according to importance sampling with uniform distribution $g(\tilde{W}_{t,t'+1})$,\STATE \hspace{13cm}  $\forall \tilde{\omega}\in\tilde{\Omega}_{t,t'+1}^e(\tilde{S}_{tt'}^{x})$
\STATE \hspace{0.5cm} $\tilde{S}_{t,t'+\tau(\tilde{x}_{tt'},\tilde{W}_{t,t'+1})}=\tilde{S}^{M,x}(\tilde{S}_{tt'}^{x},\tilde{W}_{t,t'+1})$
\STATE \hspace{0.5cm} $t'\leftarrow t'+\tau(\tilde{x}_{tt'},\tilde{W}_{t,t'+1})$
\STATE \hspace{0.2cm} \textbf{end if}
\STATE \textbf{end while}
\STATE
\STATE \textbf{function} $SimPolicy(\tilde{S}_{tt'})$
\STATE \hspace{0.5cm} Choose a sample path $\tilde{\omega}\in \tilde{\Omega}_{tt'}$
\STATE \hspace{0.5cm}\textbf{while} $\tilde{S}_{tt'}$ is non-terminal
\end{algorithmic}
\end{algorithm}

\begin{algorithm}[t!] \footnotesize
\begin{algorithmic}[t!]
\STATE \hspace{1cm}Choose $\tilde{x}_{tt'}\leftarrow \pi_0({\tilde{S}_{tt'}})$
\STATE \hspace{1cm} ${\tilde{S}_{t,t'+\tau(\tilde{x}_{tt'}(\tilde{\omega}))}}\leftarrow \tilde{S}^M(\tilde{S}_{tt'},\tilde{x}_{tt'}(\tilde{\omega}))$
\STATE \hspace{1cm} $t'\leftarrow t'+\tau(\tilde{x}_{tt'}(\tilde{\omega}))$
\STATE \hspace{0.5cm}\textbf{end while}
\STATE \textbf{return }$\tilde{V}_{tt'}(\tilde{S}_{tt'})$ (Value function of $\tilde{S}_{tt'}$)
\STATE \textbf{function} $Backup(\tilde{S}_{tt'},\tilde{V}_{tt'}(\tilde{S}_{tt'}))$
\STATE \hspace{0.5cm}\textbf{while} $\tilde{S}_{tt'}$ is not null\textbf{ do}
\STATE \hspace{1cm} $N(\tilde{S}_{tt'})\leftarrow N(\tilde{S}_{tt'}) + 1$
\STATE \hspace{1cm} $t^*\leftarrow$ time when the truck was at predecessor node, i.e., $(\tilde{S}_{tt'}=\tilde{S}^{M,x}(\tilde{S}_{tt^*}^x,{\tilde{W}_{t,t^*+1}}))$ where \STATE \hspace{12cm}$(\tilde{S}_{tt^*}^x=\tilde{S}^{M}(\tilde{S}_{tt^*},{\tilde{x}_{tt^*}}))$
\STATE \hspace{1cm} $\tilde{S}_{tt^*}^{x}\leftarrow \mathrm{~predecessor~of~} \tilde{S}_{tt'}$
\STATE \hspace{1cm} $N(\tilde{S}_{tt^*},{\tilde{x}_{tt^*}})\leftarrow N(\tilde{S}_{tt^*},{\tilde{x}_{tt^*}}) + 1$
\STATE \hspace{1cm} $\tilde{V}^x_{tt^*}(\tilde{S}^x_{tt^*})\leftarrow \frac{1}{\sum_{\tilde{W}_{t,t^*+1}\in\tilde{\Omega}_{t,t^*+1}^e(\tilde{S}_{tt^*}^{x})} p(\tilde{W}_{t,t^*+1})} \cdot E_{g}[p(\tilde{W}_{t,t^*+1})/g(\tilde{W}_{t,t^*+1})\tilde{V}_{tt^*}(\tilde{S}^{M,x}(\tilde{S}_{tt^*}^x,{\tilde{W}_{t,t^*+1}}))]$
\STATE \hspace{1cm} $\tilde{S}_{tt^*}\leftarrow \mathrm{~predecessor~of~} \tilde{S}_{tt^*}^{x}$
\STATE \hspace{1cm} $\Delta \leftarrow \tilde{C}(\tilde{S}_{tt^*},{\tilde{x}_{tt^*}}) + \tilde{V}^x_{tt^*}(\tilde{S}^x_{tt^*}) $
\STATE  \hspace{1cm}  $\tilde{V}_{tt^*}(\tilde{S}_{tt^*})\leftarrow \tilde{V}_{tt^*}(\tilde{S}_{tt^*})+\frac{\Delta- \tilde{V}_{tt^*}(\tilde{S}_{tt^*})}{N(\tilde{S}_{tt^*})}$
\STATE \hspace{1cm} $t'\leftarrow t^*$
\STATE \hspace{0.5cm}\textbf{end while}
\end{algorithmic}
\end{algorithm}

\subsection{The steps of MCTS}
Below we describe in more detail each of the four steps of MCTS.

\begin{itemize}
\item[1. ]\textbf{Selection:} In sampled MCTS, there are two selection strategies which are applied based on the domain of selection; one is  for the decision space while the other is for the exogenous event space. Starting from the root node, selection chooses  a decision based on previous gained information while controlling a  balance between exploration and exploitation. The most popular method used in the computer science literature is Upper Confidence Bounding applied to Trees (UCT)~(\citet{BPW12,KS06}). However, upper confidence bounds are used for maximization problems. In this work, the aim is to minimize the objective function which is equivalent to maximizing its negative value. UCT builds on an extensive literature in computer science on upper confidence bounding (UCB) policies for multiarmed bandit problems~(\citet{ABF02}). UCT selects the decision that maximizes the following equation:
    \begin{equation}
    \tilde{x}_{tt'}^*=\arg\max_{\tilde{x}_{tt'}\in \tilde{\mathcal{X}}_{tt'}^e(\tilde{S}_{tt'})}\left(-\left(\tilde{C}(\tilde{S}_{tt'},\tilde{x}_{tt'})+\tilde{V}^x_{tt'}(\tilde{S}^x_{tt'})\right)+\alpha\sqrt{\frac{\ln N(\tilde{S}_{tt'})}{N(\tilde{S}_{tt'},\tilde{x}_{tt'})}}\right),\label{eq:UCB}
    \end{equation}
    where $\alpha$ is a hyper-parameter that balances exploration and exploitation. The choice of decision $\tilde{x}_{tt'}^*$ that maximizes the UCT equation depends on a weighted average of two terms; the first term of the UCT equation represents the average value of the state-decision after $N(\tilde{S}_{tt'},\tilde{x}_{tt'})$ iterations. So, the higher the average value of the state-decision, the more it contributes to exploiting the decision further since its reward is high. The second term gives a higher weight to the decision that has been less explored since its value decreases as $N(\tilde{S}_{tt'},\tilde{x}_{tt'})$ increases which contributes to exploring the decisions with lower number of visits.

 \indent Upon choosing $\tilde{x}_{tt'}^*$, the state of the network becomes $\tilde{S}^{x}_{tt'}$ after which we sample an exogenous realization $\tilde{W}_{t,t'+1}$ from the set of explored exogenous events $\tilde{\Omega}^e_{t,t'+1}(\tilde{S}_{tt'})$ for state $\tilde{S}_{tt'}$.
     In our model, the exogenous events may have very different probability density functions where some events can lie on the tail of the probability density function. Thus, in order to avoid too many iterations to catch the rare events, we propose importance sampling to choose $\tilde{W}_{tt'}^*$ from the set of available samples. Importance sampling yields the same expected value of the outcome of a random variable with a much lower number of iterations compared to sampling using the random variable's initial probability density function.

      Assume that $\tilde{\omega}$ represents an outcome of the exogenous random variable $\tilde{W}_{tt'}$. Since there are only a few outcomes for $\tilde{W}_{tt'}$, let $p(\tilde{W}_{tt'}=\tilde{\omega})$ be the probability mass function for outcome $\tilde{\omega}$ and $E_{p}[\tilde{W}_{tt'}]$ be its expected value. Also, define a new probability mass function $g(\tilde{W}_{tt'}'=\tilde{\omega})$ which is designed to balance the selection of all outcome events. According to importance sampling, $E_{p}[\tilde{W}_{tt'}]\approx E_{g}[p(\tilde{W}_{tt'}')/g(\tilde{W}_{tt'}')\tilde{W}_{tt'}']$ for a much lower number of realizations of $\tilde{W}_{tt'}'$. For simplicity, let $\tilde{\omega}_{tt'}$ be an abbreviation of the event $\tilde{W}_{tt'}=\tilde{\omega}$. In this work, we choose $g(\tilde{W}_{tt'})$ to have a uniform distribution for all random events that can take place from state $\tilde{S}_{t,t'-1}^{x}$. Then, one of the outcomes $\tilde{\omega}$ is chosen according to $g(\tilde{W}_{tt'})$ and later in the backpropagation step,  the value function of $\tilde{S}^{M,x}(\tilde{S}_{t,t'-1}^{x},\tilde{\omega}_{tt'}^*)$ is weighted by $p(\tilde{\omega}_{tt'}^*)/g(\tilde{\omega}_{tt'}^*)$ in order to maintain the same expected value of the exogenous events.

\item[2. ]\textbf{Expansion:} This is the process of adding a child node  to the tree to expand it. Upon visiting a state, one can either expand an unexplored state (via a decision or exogenous event) or exploit existing states. For example, one can set a threshold, $d^{thr}$, for the number of decisions and another threshold, $e^{thr}$, for the number of exogenous events to be expanded first before starting the exploitation process.

    If a state has several unexplored decisions and exogenous events represented by the sets $\tilde{\mathcal{X}}_{tt'}^u(\tilde{S}_{tt'})$ and $\tilde{\Omega}_{tt'}^u(\tilde{S}_{tt'}^{x})$, respectively,  then an unexplored decision is chosen optimistically by using a two-stage lookahead model. That is, for each unexplored decision $x_{tt'}\in \tilde{\mathcal{X}}_{tt'}^u(\tilde{S}_{tt'})$, we generate a random exogenous event $\tilde{\omega}\in \tilde{\Omega}_{t,t'+1}^u(\tilde{S}_{tt'}^{x})$ and evaluate the value of the obtained state $\tilde{S}^{M,x}(\tilde{S}_{tt'}^{x},\tilde{\omega}_{t,t'+1})$ by calling the simulation policy. Then, the decision and its corresponding exogenous event that corresponded to the highest obtained value are chosen to be expanded in the tree. Finally, if the selection phase reaches a  state, say  $\tilde{S}_{tt'}^{x}$, that is already part of the tree but for which the threshold value of the number of exogenous events to be explores is not met, then a random exogenous sample $\tilde{\omega}\in \tilde{\Omega}_{t,t'+1}^u(\tilde{S}_{tt'}^{x})$ is created resulting in state  $\tilde{S}^{M,x}(\tilde{S}_{tt'}^{x},\tilde{\omega}_{t,t'+1})$ and evaluated as discussed above.

\item[3. ]\textbf{Simulation:} The simulation policy is typically a heuristic to provide an initial estimate of the value of the state that has just been added to the tree. Adding a node to the search tree at time $t'$, results in state $\tilde{S}_{tt'}$ and updates the probability space $\tilde{\Omega}_{tt'}(\tilde{S}_{tt'})$. The simulation policy selects states from the newly added state until the end of the simulation; this is a roll-out simulation starting from the expanded state.  The proposed simulation policy is based on information relaxation by generating a sample path $\tilde{\omega}\in\tilde{\Omega}_{tt'}(\tilde{S}_{tt'})$ to evaluate the value of $\tilde{S}_{tt'}$. The sample path determines the set of power lines  that have faulted along with the fault types, and required travel and repair times for each arc in the graph.

    {\color{black} We design a simulation policy based on information relaxation. The essence of information relaxation is to relax non-anticipativity constraints (i.e., allow the decision maker to use future information) which produces an optimistic estimate of the value of a node.
    If the value of the unexpanded node is better than what we have in the tree, then we expand it.
    }

    We elaborate on the simulation policy in Appendix~\ref{sec:sim_policy} where we formulate the problem as a sequential optimization problem; the resulting problem is a mixed integer non-linear program (MINLP) where the nonlinearity arises from the radial structure of the grid. In general, MINLP problems are difficult to solve for large network sizes but Appendix~\ref{sec:sim_policy} shows that the problem reduces to  a travelling salesman problem (TSP) since the objective is to find the tour of the truck that visits each location with a fault exactly once to repair it in order to minimize the customer outage-minutes. So, Appendix~\ref{sec:sim_policy} shows that the optimal route of the truck given a sample path can be optimally solved via dynamic programming for a small number of generated faults whereas a heuristic TSP solution becomes necessary if the number of generated faults is large.

\item[4. ]\textbf{Backpropagation:} At the end of the simulation, a value, $\tilde{V}_{tt'}(\tilde{S}_{tt'})$, of the newly created state, $\tilde{S}_{tt'}$, is obtained. Starting from the last added node in the tree, its simulated value is backpropagated through all ancestors of state $\tilde{S}_{tt'}$  until the root state to update their statistics. The number of visit counts of all  ancestor states of state $\tilde{S}_{tt'}$ are increased by one and their values are modified according to a chosen criteria where we choose the average value of all rollouts through a state. While backpropagating, assume that we have $\tilde{S}_{tt'}=\tilde{S}^{M,x}(\tilde{S}_{tt^*}^x,{\tilde{\omega}_{t,t^*+1}})$, then the value of the ancestor post-decision state, say $\tilde{S}_{tt^*}^x$, should be updated as $\tilde{V}^x_{tt^*}(\tilde{S}^x_{tt^*})\leftarrow \frac{1}{\sum_{\tilde{W}_{t,t^*+1}\in\tilde{\Omega}_{t,t^*+1}^e(\tilde{S}_{tt^*}^{x})} p(\tilde{W}_{t,t^*+1})} \cdot E_{g}[p(\tilde{W}_{t,t^*+1})/g(\tilde{W}_{t,t^*+1})\tilde{V}_{tt^*}(\tilde{S}^{M,x}(\tilde{S}_{tt^*}^x,{\tilde{W}_{t,t^*+1}}))]$ where the expectation is over all the explored exogenous events from post-decision state $\tilde{S}_{tt^*}^x$. Also, the value function of the ancestor pre-decision state, $\tilde{S}_{tt^*}$ should be updated with a value of $\Delta \leftarrow \tilde{C}(\tilde{S}_{tt^*},{\tilde{x}_{tt^*}}) + \tilde{V}^x_{tt^*}(\tilde{S}^x_{tt^*})$ which accounts for the link cost and the updated post-decision state value function so that we get
  $\tilde{V}_{tt^*}(\tilde{S}_{tt^*})\leftarrow \tilde{V}_{tt^*}(\tilde{S}_{tt^*})+\frac{\Delta- \tilde{V}_{tt^*}(\tilde{S}_{tt^*})}{N(\tilde{S}_{tt^*})}$ (note that we assume that the weight of all decisions from the same state is equal).
\end{itemize}

\section{Experimental Results}\label{sec:results}
To assess the performance of the proposed approaches, the simulated power grid is constructed using real data provided by PSE\&G which describes the structure of circuits in their electrical distribution network. The data corresponds to the northeastern portion of PSE\&G's power grid in New Jersey and it is formed of $319$ circuits, each rooted at a substation that connects the circuit (that is, the distribution grid) to the high voltage transmission grid. There are an average of $41$ protective devices and $724$ power lines per circuit. The data identifies the type and location of each component in the circuits such as substations, protective devices, power lines and transformers.

The simulator is programmed to generate storms that pass across the grid generating power line faults causing total or partial circuit power outages~(\citet{ABP14}). {\color{black} All simulations are assumed to last 24 hours, and the lookahead policy always plans to the end of the horizon.} The obtained outages trigger some of the affected customers to call to report the outage. When a power line faults due to storm damage, the simulator finds the nearest upstream protective device, opens it (shutting off power), and then identifies all the customers who subsequently lose power. For each customer experiencing a power outage, a Bernoulli random variable is generated, with a probability of success which equals the call-in probability across the customer's segment,  to determine whether the customer will call or not. Thus, the higher the call-in probability, the higher the number of lights-out calls.

Recall that to repair the grid due to storm damage, a truck uses  a roadway defined by the graph $G(\mathcal{V},\mathcal{E})$ where $\mathcal{E}$ is the set containing the arcs/roads between two consecutive nodes of the graph. 
We aggregate the power network into segments that group lines and poles that are served by the same protective device.  The truck is routed from one segment to the next over the road network, where the truck is assumed to stop and repair any damage within the same segment at the same time, since they are likely to be easily visible once the truck is in the area.

Basically, the network  $G(\mathcal{V},\mathcal{E})$ is formed of the poles that carry the protective devices where $\mathcal{E}$ corresponds to the connection matrix in the form of the minimum distance between any two nodes in the network according to the real roadway.

As the storm passes across the grid, each power line $i$ along its way is associated with a prior probability of power line fault as explained in detail in~\citet{ABP14}. In the simulator, the prior probability of fault of the power lines are generated based on several parameters such as the severity of the storm, its diameter and the distance of the power line from its center. Since this paper addresses a single truck, we tuned the priors to create storms that generated one to as many as tens of outages. In this case, the segment fault probabilities range between  $0$ and $0.765$ where the segments that faulted have posteriors ranging between $0.032$ and $0.765$.

\subsection{The Lookahead Policy}
After collecting the priors for power line faults and the customer calls, the simulator executes the proposed power line fault probability model presented in (\ref{eq:eq3}) upon which the utility truck is routed to restore the grid. According to the simulations, the average number of segments affected by the storm path is $1558$. Based on the statistics of $1000$ networks, the minimum posterior fault of a segment that faulted is $0.032$, so we choose the segments with posterior probability of faults greater or equal to $0.01$ as candidate segments that have faulted.

\begin{figure}[t!]
\centering
\includegraphics[scale=0.56]{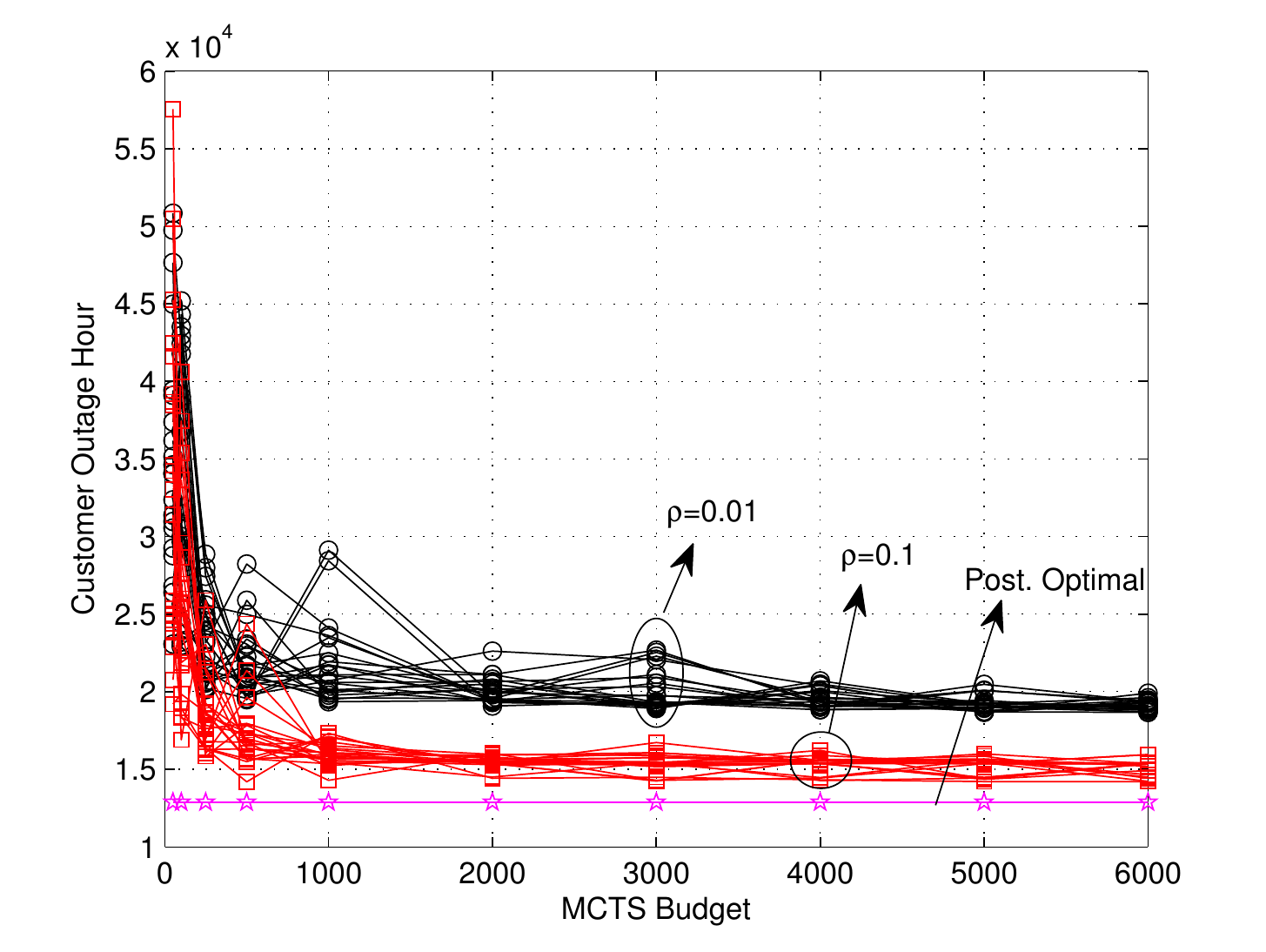}
\caption{\label{fig:Simu_COH} Customer outage-hours vs. MCTS budget for one network; twenty simulations for various call-in probabilities ($\rho$) where each simulation corresponds to the route of the truck starting from the depot until the stopping condition is met.}
\end{figure}

Figures~\ref{fig:Simu_COH} shows the customer outage-hours vs. MCTS budget for twenty simulations for one chosen network which has five faults. {\color{black} The graphs are shown for each simulation to indicate the level of variability.}  The MCTS budget refers to the number of iterations $n^{thr}$ executed to compute the MCTS tree, while a simulation means running the truck until the stopping criterion in Algorithm~\ref{alg:adp} is met. 


Since MCTS depends on Monte Carlo sampling, it becomes more consistent as the computation budget increases. Figure~\ref{fig:Simu_COH} illustrates the behavior of the policy as a function of the MCTS budget, for two values of the customer call-in probability rho:  $\rho=.01$, and $\rho = 0.10$.  As more customers call in, we learn more about the network which reduces the uncertainty. The results demonstrate that if $\rho=0.10$, an MCTS budget of 1000 provides consistently reliable results; if $\rho=0.01$, {\color{black} we need a budget of 2000 to get consistently high quality results, with sharply diminished returns with larger budgets.} The results are also compared to the posterior optimal solution which corresponds to the optimal solution after revealing the locations of faults in the network  obtained using the dynamic program presented in Algorithm~\ref{alg:dptsp}.


{\color{black} We tune the hyper-parameter $\alpha$ (presented in (\ref{eq:UCB})) by performing a brute force search for the best value of $\alpha$ over the range from $0.1$ till $7$ with increments of $0.1$. The minimum objective is obtained for $\alpha=2.2$. So, for the rest of the simulations, we set the hyper-parameter $\alpha = 2.2$.}

\begin{figure}[t!]
\centering
\includegraphics[scale=0.55]{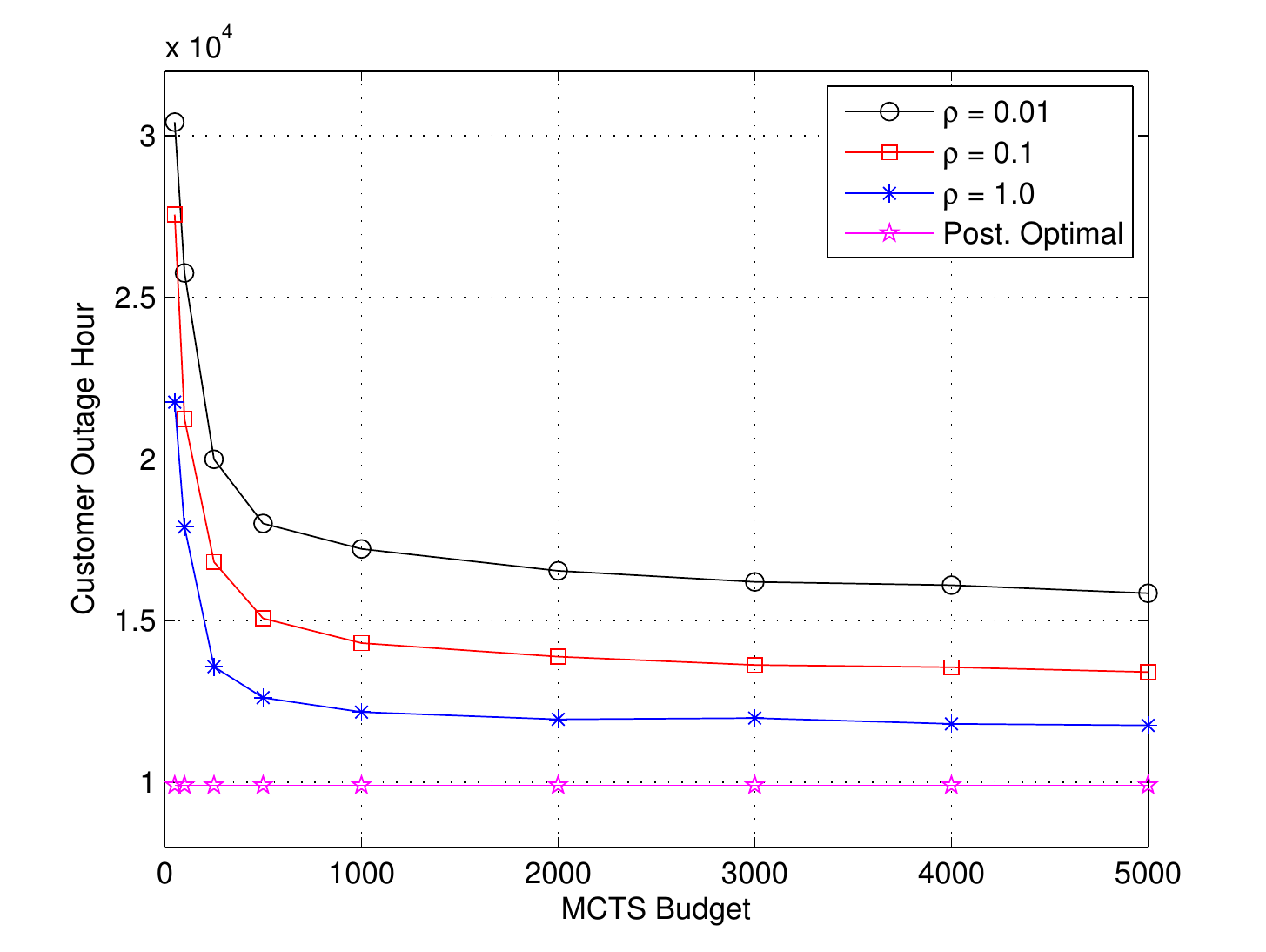}
\caption{\label{fig:Av_COH}Average customer outage-hours vs. MCTS budget for ten networks; for each network, twenty different storms are simulated and the resulting customer outage-hours of the $200$ simulations are averaged.}
\end{figure}

Figure~\ref{fig:Av_COH} shows the average customer outage-hours vs. MCTS budget of ten networks for various call-in probabilities; for each network, twenty different storms are simulated and the resulting customer outage-hours of the $200$ simulations are averaged.  {\color{black} We only show the averages, but the sample paths in figure~\ref{fig:Simu_COH} indicate that the level of variation is relatively low, and the averaging over 200 simulations produces an accurate estimate of the mean.} The number of faults in the networks ranges from $4$ to $12$ with an average of $6.09$ faults.

From these experiments, we make the following observations:
\begin{itemize}
\item As the  call-in probability increases, the lookahead policy provides a solution that is closer to the posterior optimal since more information is provided. Also, the required MCTS budget decreases as the call-in probability increases; that is, for a call-in probability of $0.01$ and $1.0$, around $4000$ and $1000$ MCTS iterations are required, respectively, to converge to a good solution. The lookahead policy provides a solution that is $18.7$\%, $35.3$\% and  $58.5$\% higher than the posterior optimal for a call-in probability of $1.0$, $0.1$ and $0.01$ respectively.  {\color{black} Note that asymptotic optimality of optimistic MCTS does not mean that these gaps will go to zero.  There will always be a gap as long as the underlying problem is stochastic.}
\item It is revealed that a high gain is obtained when the call-in probability increases from $0.01$ to $0.1$ since the information provided can be from different locations which help in detecting the location of outages. As the call-in probability rises above $0.1$, the benefits are more modest than the increase from $0.01$ to $0.1$ since we only need one customer out of a group across the same segment to make the phone call.
\end{itemize}
\begin{figure}[t!]
\centering
\includegraphics[scale=0.55]{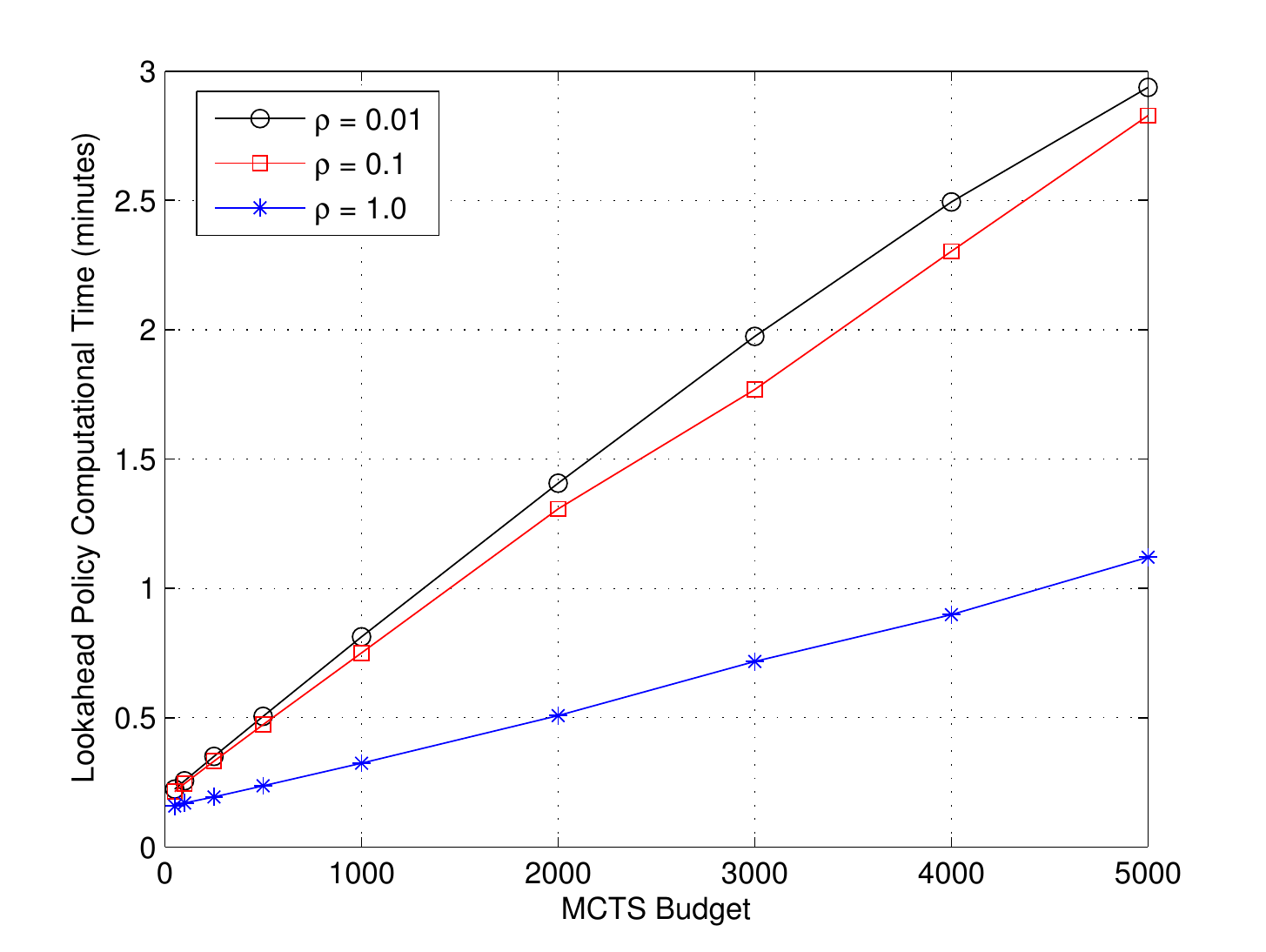}
\caption{\label{fig:Av_CT}Average computational time of the lookahead policy to decide on the truck's next hop vs. MCTS budget.}\vspace{-0.4cm}
\end{figure}

Figure~\ref{fig:Av_CT} shows the average computational time of the lookahead policy to decide on the truck's next move. Obviously, the computational time increases as the MCTS budget increases; however,  no additional gain in terms of the objective value is obtained.  The computational time decreases as the call-in probability increases since more information is provided. Figure~\ref{fig:Av_COH} shows that the required  computational budget is $4000$ for a good performance; this corresponds to a computational time of $2.5$ minutes which is suitable for online problems.

{\color{black}
Caution needs to be used when evaluating computational requirements.  This work has to be viewed as the first stage in the introduction of a new algorithmic technology for a new problem.  CPU times can be reduced, possibly significantly, through strategies such as warm starts (using the solution from $t-1$ to accelerate the solution for time $t$), and parallelism (different paths in the tree can be explored simultaneously).  In addition, we need to consider the rich array of strategies for approximating stochastic lookahead models to find the best strategy for this particular application.
}

\begin{figure}[t!]
\centering
\includegraphics[scale=0.55]{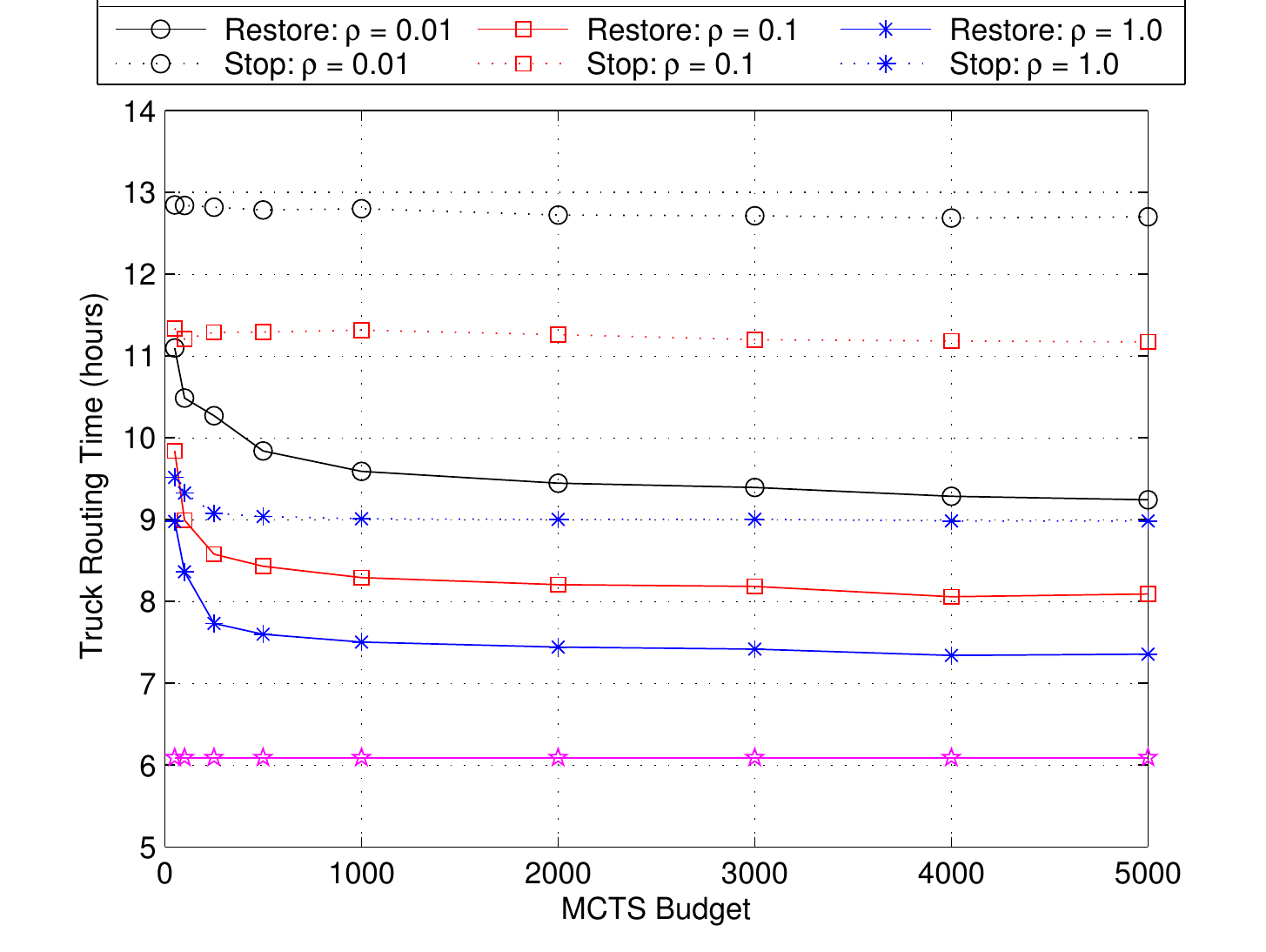}
\caption{\label{fig:AV_TRT} Average truck routing time to restore the grid and to stop routing vs. MCTS budget; the truck is stopped when the stopping condition of the lookahead policy is met.}\vspace{-0.0cm}
\end{figure}

Figure~\ref{fig:AV_TRT} shows the time required to restore the grid and the time required to stop truck routing compared to the posterior optimal solution.  Even when the grid is completely restored,  the utility cannot  detect that unless the distribution system is fully equipped with sensors which is not the case. In this work, the utility center relies on the probability model to decide when the grid is restored and consequently, stop the truck routing process.  The posterior optimal solution indicates that  $6$ hours are required to restore the grid if full information is revealed. For a call-in probability of $1.0$, the lookahead policy can restore the grid in $7.2$ hours on average, but the truck is stopped after $11$ hours. For a lower call-in probability, both the restore and stop times increase for the same reasons mentioned previously but still have a good performance with respect to the posterior optimal~solution.

%

\subsection{Industrial Heuristics}
We also compare the performance of the proposed  lookahead policy to industrial heuristics which typically rely on escalation algorithms based on our discussion with PSE\&G. An escalation algorithm is a simple policy function approximation formed of predefined decision rules. The basic idea is that it back-traces from each lights-out call location to find the first common point for all the calls.

Escalation algorithms are good for locating a single fault which is assumed to be an upstream fault that triggers all downstream calls. However, this is often not the case as there can be more than one fault triggering the calls.  Escalation is performed at the control center along with other intelligence techniques as explained in the literature review in  Section~\ref{sec:lit_rev}.

For a fair comparison, we have constructed an industrial heuristic based on escalation that relies on the information that we are using which includes the lights-out calls, the grid structure, and the storm path. Thus, the proposed escalation algorithm gives priority to visit the common locations that would trigger all calls but after that searches for other downstream faults as explained in Algorithm~\ref{alg:esc}.

\begin{algorithm}[t!] \footnotesize
\caption{Escalation Algorithm for Grid Restoration}
\label{alg:esc}
\begin{algorithmic}[t!]
\STATE\textbf{Step 1.} \textbf{For }each circuit \textbf{do}
\STATE\hspace{1.1cm}\textbf{Step 1a.} Collect all calls and back trace to find the first node that is common to all calls say node $x$.
\STATE\hspace{1.1cm}\textbf{Step 1b.} Send the truck to node $x$ and then back trace to the substation to cover all upstream faults \STATE\hspace{1.5cm}(when a truck visits a node, it fixes an existing fault and this applies to all steps of the algorithm).
\STATE\hspace{1.1cm} \textbf{Step 1c.} From node $x$, perform down tracing to reach the first segment from which a call was initiated \STATE\hspace{2.5cm}  and place it in set~$\mathcal{D}$.
\STATE \textbf{Step 2.} \textbf{For} each segment in $\mathcal{D}$ \textbf{do}
\STATE\hspace{1.1cm} \textbf{Step 2a.} Perform down tracing to cover all nodes that called.
\end{algorithmic}\vspace{-0.0cm}
\end{algorithm}

\begin{table}[t]
\vspace{-0.0cm}
\caption{Average Statistics for Escalation Algorithm}\label{tab:esc}
\begin{center}\vspace{-0cm}
\begin{tabular}{|l|c|c|c|}
\hline
 & $\rho = 0.01$ & $\rho = 0.1$ & $\rho = 1.0$   \\
\hline
Customer outage-hours & $ 2.34*10^4$& $2.21*10^4$ & $2.17*10^4$\\
\hline
Number of unrepaired faults& $1$ &$0.42$ &$0$ \\
\hline
Number of customers in outage & $14$ & $ 5.9$ & $0$ \\
\hline
Time to restore (hours)& $21.83$ &  $20.76$ & $20.34$ \\
\hline
Time to stop
 (hours) & $48$ & $48$ & $26.46$ \\
\hline
\end{tabular}
\end{center}\vspace{-0.0cm}
\end{table}

Table~\ref{tab:esc} shows the average statistics of the same ten networks used earlier using $200$ simulations. In the escalation algorithm, searching for a fault is triggered by the customer calls because there is no clue for the utility center to predict the locations of faults except by visiting the locations that  trigger the calls if faulted. So, if there is a segment with a low number of customers where no one called to report an outage, then there is no way to know if a fault is present.  By contrast, our probability model may still return a positive probability of an outage.

In Table~\ref{tab:esc}, we see that for a call-in probability of $0.01$, on average one fault could not be identified since a total of $14$ customers are only affected and none of them called at such a low call-in probability. So, the utility center will stop routing the truck assuming that it has recovered the grid, because it visited every location from which a call was received. The average number of faults is $6.1$ with an average of one unrepaired fault. So, the time it to took the truck to restore an average of $5.1$ faults that triggered calls is around $21.83$ hours and the customer outage-hours is $2.34*10^4$.

We set the maximum truck routing time to $48$ hours unless it restored the entire grid in a smaller amount of time. The escalation algorithm was fast in locating an upstream fault whenever there is one that triggered all downstream calls; however, it took a long time to locate downstream faults because there is no clue from the common point where a fault might be except by following the locations that could trigger the calls. The lookahead policy, on the other hand, is able to repair all faults in 9.15 hours  with customer outage-hours of $1.58*10^4$ showing its superior performance compared to the industrial heuristic.

The lookahead policy outperforms the industrial escalation policy both through more efficient learning, and the ability to explore streets even when there is no phone call that would indicate that an outage may have occurred there.  The restore times drop from 20.3 to 7.2 hours using MCTS, while the time when the vehicle stops looking drops from 26.5 to 9.0 hours.  Total customer outage-hours drops from $2.17*10^4$ to $1.17*10^4$ using MCTS.

{\color{black}
\section{Areas for further research}\label{sec:conclusion}


This paper has proposed a new class of routing problem that involves the management of a vehicle that simultaneously performs physical tasks (such as repairing segments of the grid, but other tasks could be substituted) while simultaneously learning uncertain parameters.  The result is a dynamic system with both a physical state and a belief state.

The problem of managing a resource that can perform both learning and mitigation represents a new class of routing problem that involves both a physical state and a belief state (possibly along with other information).  We build on a general way of modeling sequential decision problems, highlight all four classes of policies (from \citet{Powell2022RLSO}), and then propose to use a stochastic lookahead policy, which is easily the hardest of the four classes.  We then propose a practical, scalable strategy based on Monte Carlo tree search, which naturally handles high-dimensional state variables.

Stochastic lookahead policies based on equation \eqref{eq:opt_policy} represent an important area of research in transportation and logistics, since these policies are fundamentally intractable.  There are two broad lines of investigation: a) choosing tractable classes of policies for the ``policy within a policy,'' and b) replacing the lookahead model with a simpler approximation.

Some research ideas for policies for the ``policy within a policy'' of equation \eqref{eq:opt_policy} include:
\begin{itemize}
\item Specific applications will introduce special structure that will suggest simple rules (PFAs in the language of section \ref{sec:designingpolicies}) that will be easy to implement, fast to compute and yet produces good results.  These all fall under the heading that we have called ``pessimistic MCTS,'' and it would be interesting to compare these to other approaches such as the ``optimistic MCTS'' used in this paper.
\item Approximate dynamic programming can be useful for certain classes of problems, where we use VFA-based policies.
\item Deterministic lookaheads (DLAs) which can be based on forecasts, or use the concept of information relaxation.
\item Hybrid policies such as the ``offline-online'' policy of \cite{Ulmer2019} which uses a VFA-based rollout policy to produce online estimates of the value of being in a state.
\item We need theoretical insights into the behavior of different classes of policies.
\end{itemize}

We anticipate that the richest areas of research will lie in the replacement of the complete base model with an approximation that is easier to compute yet works well.  It is clear that this research will be highly problem dependent and will need to be conducted in the context of major problem classes.  Some research directions in this area include:
\begin{itemize}
\item Compare modeling the full information process (such as the lights-out calls as well as observations of the outages) against more limited processes (such as just the outages).
\item Compare using a fixed belief state (which eliminates learning in the lookahead model) with a dynamic belief state, using any of the information processes above.  This would produce comparisons of active learning (updating the belief state) against passive learning (if the belief state is held constant).
\item There will be settings (such as ours) where updating beliefs require the use of expensive applications of Bayes' theorem.  Investigate simpler updating to streamline this expensive calculation.
\item As with policies, we need research that provide insights into the effect of different types of approximations (e.g. deterministic lookahead) on the performance of lookahead policies.
\end{itemize}

Finally, there are some important problem extensions worth considering:
\begin{itemize}
\item As of this writing, MCTS is limited to managing a single resource at a time (the decision $x_t$ cannot be a vector), but there are many settings where we have to manage fleets of vehicles, people or devices.
\item Related to the first extension, an important problem class is multiagent systems, where different agents perform these functions independently, but with information sharing.  This opens up the door for different forms of information sharing to achieve coordination, without having to solve a high-dimensional resource allocation problem.
\end{itemize}
}

\setlength{\bibsep}{1.0pt}

\clearpage
\section{Appendices}
\subsection{Components of the Distribution Power Grid}\label{sec:DPG}
We provide a brief overview about the components of the distribution power grid for readers who are not familiar with this background.

The \textbf{distribution grid} is the final stage in the delivery of electric power. It carries electricity from the transmission system to individual consumers. It mainly has a radial structure.  Distribution substations connect to the transmission system and lower the transmission voltage to medium voltage with the use of \textbf{transformers}. Primary distribution \textbf{power lines} carry this medium voltage power to distribution transformers located near the customer's premises. Distribution transformers again lower the voltage to the utilization voltage used by appliances.

An \textbf{electric circuit} is the set of electric components and power lines in which electricity flows from a source (power generating station) to a destination e.g., consumers. An overhead distribution grid is carried on poles and each pole can carry several independent circuits.

Power systems contain \textbf{protective devices} to prevent injury or damage during failures. The quintessential protective device is the fuse. When the current through a fuse exceeds a certain threshold, the fuse element melts, producing an arc across the resulting gap that is then extinguished, interrupting the circuit. Given that fuses can be built as the weak point of a system, fuses are ideal for protecting circuitry from damage.

{\color{black}
Using the configuration of the grid, we define node $d_{ui}$ to be the first upstream protective device of node $i$  on circuit~$u$. For each circuit~$u$, let $\mathcal{D}_u$ be the set of protective nodes  and $\mathcal{S}_u$ be the set of segments where each segment contains the power lines that trigger the same protective device. We define ${\cal S}_{ui}$ as the set of power lines belonging to the same segment of node~$i$ on circuit~$u$, i.e., triggering the same protective device, $\mathcal{S}_{ui}=\{j, d_{uj} = d_{ui} \}$. We also introduce $\mathcal{Q}_{ui}=\{j, \mbox{ node } i\in \mathcal{I}_u \hbox{~becomes in outage if power line~} j \in \mathfrak{L}_u \hbox{~faults}\}$ as the set of power lines that if faulted result in an outage to node $i$ on circuit~$u$ and $\mathcal{W}_{ui}$ as the set of downstream segments of node~$i$ but with different upstream protective devices, i.e., $\mathcal{W}_{ui}=\{S_{uj}, \hbox{~segment~}S_{uj} \hbox{~is downstream of} \linebreak  \hbox{node~}i ~\hbox{~\&~}  d_{ui}\neq d_{uj}\}$. We assume that $n_{ui}$ customers are attached to node $i$ on circuit~$u$. If the node is a transformer, then $n_{ui}>0$, otherwise $n_{ui}=0$ because no customers are attached to the power generator or protective devices.
}

\subsection{Lookahead Simulation Policy}\label{sec:sim_policy}
In the lookahead simulation policy, we evaluate the value of the newly added node at time $t'$ by using a lookahead model in which a sample path $\tilde{\omega}\in\tilde{\Omega}_{t,t'+1}(\tilde{S}_{tt'})$  is generated. The sample path determines the set of power lines  that have faulted along with the fault types, required repair and travel times for each arc in the graph.  For example, using the probability of fault of each power line in the system, we generate a Bernoulli random variable with a probability of success equals to the probability of fault. In this case, each power line in the power system has a posterior probability of fault equal to either $1$ or $0$. Thus, we define the following indicator function:
\begin{eqnarray}
\mathds{1}_{\tilde{L}_{t'ui}}=\left\{
\begin{array}{c l}
    1 & \mbox{,~if~} \tilde{L}_{tt'ui}(\tilde{\omega})=1,\\
    0 & \mbox{, otherwise},
\end{array}\right.
\end{eqnarray}
where $\mathds{1}_{\tilde{L}_{t'ui}}=1$ if  power line $i$ on circuit~$u$ is in fault at time $t'$ and it is equal to zero, otherwise.

In the simulation policy, we use the time index $t''$ for each state, decision and random variable. Recall, at time~$t$, we are in the base model where we fix the set of calls and call MCTS to find the truck's next hop. For each node included in the MCTS tree, we index it with time $tt'$.
In the fourth step of MCTS, to evaluate the value of a newly added state $\tilde{S}_{tt'}$, we use another lookahead model to generate a sample path $\tilde{\omega}\in\tilde{\Omega}_{t,t'+1}(\tilde{S}_{tt'})$ at $t'$. While, routing the truck according to the generated path, we index the nodes  by $t''$; for example, $\tilde{S}_{t''}$ indicates the state at time $t''$ in the simulation step which is used to evaluate the value of the generated state in MCTS at time $t'$.

In the simulation step, if at time $t''$, the truck visits a location that was identified as a location with fault at time $t'$, its indicator function is set to $0$ from time $t''+1$ and on. Thus, for each power line that faulted  the following relation holds:
\begin{eqnarray}
\mathds{1}_{\tilde{L}_{t''uj}}=1-\sum_{i\in\mathcal{V}}\sum_{\hat{t}=t'}^{t''-1}\tilde{x}_{\hat{t}ij}, \mbox{~if~} \tilde{L}_{tt'uj}(\tilde{\omega} )=1.\label{eq:const5ex}
\end{eqnarray}

Whereas, if power line $i$ on circuit~$u$ did not fault in the considered scenario then, $\mathds{1}_{\tilde{L}_{t''ui}}=0$ , $\forall t''\geq t'$.
Given, a sample path $\tilde{\omega}$, the aim is to find the optimal truck's route that minimizes the customer outage-minutes. Let $\tilde{C}_{t''j}$ be a random variable representing the number of customers that regain power by visiting node $j$ at time $t''$ according to scenario $\tilde{\omega}$. Then, the value of the newly added node is obtained by solving the  following optimization problem:

\begin{eqnarray}
&&\hspace{-1.0cm} \min_{\tilde{x}_{t''}} \sum_{\hat{t}=t'}^{t'+H}\left(N-\sum_{t''=t'}^{\hat{t}}\sum_{u\in\mathcal{U}}\sum_{j=1}^N\tilde{C}_{t''uj}\right)\label{eq:objective}\\
&&\hspace{-1cm} \mbox{subject to}\nonumber\\
&&\hspace{-1cm}\tilde{C}_{t''uj}=\sum_i \left(\prod_{k\in\mathcal{Q}_{uj}\backslash j}\hspace{-0.2cm}1-\mathds{1}_{\tilde{L}_{t''uk}}\right)\mathds{1}_{\tilde{L}_{t''uj}}\left(\sum_{k\in \mathcal{S}_{uj}} n_{uk} + \sum_{s\in \mathcal{W}_{uj}}\left(\prod_{w=\min\{\mathcal{W}_{uj}\}}^s \prod_{k\in w}1-\mathds{1}_{\tilde{L}_{t''uk}}\right)\sum_{k\in s}n_{uk}\right)\tilde{x}_{t''ij}, \nonumber\\ &&\hspace{12.5cm}\forall j\in\mathcal{V}, \forall t''\\ \label{eq:const1} 
&&\hspace{-1cm} \mathds{1}_{\tilde{L}_{t''uj}}=1-\sum_i\sum_{\hat{t}=t'}^{t''-1}\tilde{x}_{\hat{t}ij},\mbox{~such~that~} \tilde{L}_{tt'uj}(\tilde{\omega})=1 ,\forall j \in \mathcal{I}_u,\forall u\in\mathcal{U}, \forall t''\label{eq:const5}\\
&&\hspace{-1cm}\tilde{\Delta}_{t''ij}\geq T_{ij}(\tilde{\omega})\tilde{x}_{t''ij} + \sum_u R_{uj}(\tilde{\omega}) \left(\tilde{x}_{t''ij} -\sum_{\hat{t}=t'}^{t''-1}\tilde{x}_{\hat{t}ji} \right),\forall (i,j)\in \mathcal{E},\forall t''\label{eq:const2}\\
&&\hspace{-1cm}\tilde{\xi}_{t''j}\geq \tilde{\xi}_{t''-1i}+ \sum_i \tilde{\Delta}_{t''ij} ,\forall (i,j)\in\mathcal{E}, \forall t''\label{eq:const3}\\
&&\hspace{-1cm}\tilde{\xi}_{t''j}\leq t''\sum_i\tilde{x}_{t''ij}+\zeta\left(1-\sum_i\tilde{x}_{t''ij}\right),\forall j\in\mathcal{V}, \forall t''\label{eq:const4}\\
&&\hspace{-1cm}\sum_{t''=t'}^{t'+H}\tilde{x}_{t''ij}\leq 1,\forall (i,j)\in \mathcal{E}\label{eq:const7}\\
&&\hspace{-1cm} \sum_k \tilde{x}_{(t''+T_{jk}(\tilde{\omega}))jk}+\sum_k \tilde{x}_{(t''+T_{jk}(\tilde{\omega})+\sum_uR_{uk}(\tilde{\omega}))jk}\leq \sum_{i}\tilde{x}_{t''ij}\leq 1,\forall j\in\mathcal{V},\forall t''\label{eq:const8}\\
&&\hspace{-1cm} \tilde{C}_{t''uj}\geq 0, \tilde{\xi}_{t''j}\geq0, \tilde{\Delta}_{t''ij}\geq0, \tilde{x}_{t''ij}\in\{0,1\} \label{eq:const11}
\end{eqnarray}

This problem is a mixed integer non-linear programming (MINLP) problem. The objective~(\ref{eq:objective}) minimizes the customer outage-minute which is equivalent to maximizing the number of customers with restored power (also referred to as served customers) up to time $t$ represented by the inner sum in the objective.
Constraint (\ref{eq:const1}) determines the  number of served customers when the truck goes from its current  location, say node $i$, to node $j$ at time $t''$. The number of served customers depends on whether there is a fault upstream to node $j$ or on its segment, i.e., in set $\mathcal{Q}_{uj}$. Note that, node $j$ will be favored to be visited if there is a fault across power line $j$, i.e., if $\mathds{1}_{\tilde{L}_{t''uj}}=1$.
Depending on the structure of the power grid, if a location faults, it causes outage to the customers attached to its segment and all downstream segments. Thus, if a fault is fixed, then all these customers will be affected. But, this also depends on whether there is a fault on any downstream location as shown in~(\ref{eq:const1}). Constraint~(\ref{eq:const1}) also reveals that the number of customers by visiting power line $j$  is positive if $\mathds{1}_{\tilde{L}_{t''uj}}=1$; however, after visiting this location, say at time $t^*$, $\mathds{1}_{\tilde{L}_{t''uj}}=0, \mbox{~for~} t''> t^*$. Thus, if the truck will come across the same location for the second time, then the gain will be $0$ which favors the truck not to visit the same location more than once unless there is no other route for it. Constraint (\ref{eq:const5}) is the same as (\ref{eq:const5ex}) and it has been explained in details above.

Constraint (\ref{eq:const2}) determines the required time to traverse arc $(i,j)\in \mathcal{E}$. If there is no power line to be investigated at node $j$, then $\mathds{1}_{\tilde{L}_{t''uj}}=0$ which means that the required traversal time is equal to the  travel time which depends on the traffic conditions only. However, if there is a power line across arc $(i,j)$, then one of two rules apply; if there is a fault on power line $j$ according to sample path $\tilde{\omega}$,  then the required traversal time accounts for the travel and repair times for power line $j$. However, if arc $(i,j)$ is traversed for the second time at time $t$ then, the traversal time is just equal to the travel time since the fault was repaired when the arc was traversed for the first time.

Constraints (\ref{eq:const3})-(\ref{eq:const4}) guarantee that the truck is at node $j$ at time $t''$, only if $\tilde{\xi}_{t''j}=t''$ which sets $\tilde{x}_{t''ij}=1$. If $\tilde{x}_{t''ij}=1$, then $\tilde{\xi}_{t''j}=t''$, otherwise $\tilde{\xi}_{t''j}$ is less than a large positive number $\zeta$ as shown in (\ref{eq:const4}) but larger than the time where the truck was lastly as indicated by~(\ref{eq:const3}). But, since the objective is maximizing the number of served customers over time, the optimization problem will set the time to the least possible value that satisfies all constraints.  In (\ref{eq:const3}), $\tilde{\xi}_{t''j}$ can be equal to $t''$ only if it satisfies the required traverse times; the required time to reach node $j$ depends on the elapsed time to reach its direct predecessor, say node $i$, in addition to the required time to traverse node $j$ from node $i$, i.e., $\tilde{\Delta}_{t''ij}$.

Note that, since the objective minimizes the customer outage-minutes, then the optimization problem keeps on routing the truck to cover all power lines that faulted, as favored by (\ref{eq:const1}), until all faults are fixed. Constraint (\ref{eq:const7}) guarantees that all arcs in the graph can be visited once in one direction and consequently, at most twice (forward and backward) which is a sufficient condition to have an Eulerian path where each power line and node with positive fault probability can be visited once. Constraint~(\ref{eq:const8}) indicates that the truck can go from node $j$ to node $k$ at time $t''+\tilde{\Delta}_{t''jk}$ only if it was at node $j$ at time $t''$ and only if the necessary traverse time, $\tilde{\Delta}_{t''jk}$, has elapsed  which depends on (\ref{eq:const2}). Moreover, this constraint removes the sub-tours in the network since the truck must have visited a node before it can travel from~it. Finally, constraint~(\ref{eq:const11}) shows that all variables are positive except $\tilde{x}_{t''ij}$ which is binary.

The formulated  optimization problem is very complex mainly because it is an MINLP problem which combines the complexity of non-linear programming and integer programming both of which lie in the class of NP-hard problems.  Thus, achieving the optimal global solution is most probably never attainable for large network sizes. While there has been a tremendous achievements in solving integer programming problems given that their continuous relaxation is convex, solving non-linear optimization problems is still a non mature area that gets stuck at local optimums.

The only constraint that cannot be linearized is (\ref{eq:const1}); it can be seen that the order of non-linearity depends on the number of faults upstream and downstream of a node which is scenario dependent. Thus,  the radial structure of the grid   is the main complicating factor in the optimization problem. 

Though the problem is non-linear, the optimal solution can be attained by using dynamic programming. For each sample path, we can transform the problem into a complete graph with nodes $\mathcal{V}^f$ which contains all the power lines that have faulted and the location of the truck indexed with $0$. The connection cost between the nodes of $\mathcal{V}^f$ are calculated by summing the shortest travel time between the nodes according to $\tilde{\omega}$. Let $S$ be the set of the nodes visited by the truck and $f(S)$ a function returning the number of customers still in outage after visiting the nodes of $S$. The aim of the problem is to find the optimal sequence of the truck route that visits each node exactly once (to repair it) in order to minimize the customer outage-minutes. This problem is equivalent to a travelling salesman problem (TSP) which is NP-Complete; however the solution can be obtained optimally using dynamic programming with complexity $O(n^22^n)$ where $n$ is the number of nodes in the TSP graph which corresponds to the number of generated faults. However, since the number of generated faults is relatively small (less than 20 faults), obtaining the optimal solution using dynamic programming is computationally  feasible.

Let $C^f(S,i)$ be the customer outage-minutes of going from vertex $0$ through the nodes of $S$ ending at node~$i$. Then, the recurrence relation of the dynamic program by going from node~$i$ to node~$j$ can be defined~as
\begin{eqnarray}
C^f(S,j) =  C^f(S-\{j\},i)+f(\{S-\{j\}\})\cdot(T_{ij}+\sum_uR_{uj}),
\end{eqnarray}
where the first term of the summation accounts for the customer outage-minutes up to node~$i$ whereas the second term accounts for the  cumulative customer outage-minutes by going from node~$i$ to node~$j$. The detailed steps of the dynamic program to obtain the value of the objective function are presented in Algorithm~\ref{alg:dptsp}.

\begin{algorithm}[t!] \footnotesize
\caption{Dynamic program for optimal customer outage-minutes of a given sample path}
\label{alg:dptsp}
\begin{algorithmic}[t!]
\STATE\textbf{Step 0.} \textbf{Initialization}
\STATE \hspace{1.2cm} \textbf{For }all $j\in\mathcal{V}^f, j\neq 0$ \textbf{do}
\STATE \hspace{1.4cm} $C^f(\{0,j\},j) = f(\{0,j\})\cdot(T_{0j}+\sum_uR_{uj})$
\STATE\textbf{Step 1.}~ \textbf{Compute customer outage-minutes for all subsets}
\STATE \hspace{1.2cm}\textbf{For} $s=3$ \textbf{to }$|\mathcal{V}^f|$
\STATE \hspace{1.2cm} \textbf{For }all subset of $\mathcal{V}^f$ of size $s$ \textbf{do}
\STATE \hspace{1.4cm} \textbf{For} all $j\in S, j\neq 0$
\STATE \hspace{1.6cm}   $C^f(S,j) = \min_{i\in S, i\neq j} C^f(S-\{j\},i)+f(\{S-\{j\}\})\cdot(T_{ij}+\sum_uR_{uj})$
\STATE \textbf{Step 2.} \textbf{Optimal solution}
\STATE \hspace{1.2cm}  $\min_{j\in\mathcal{V}^f}C^f(\mathcal{V}^f,j)$
\end{algorithmic}\vspace{-0.0cm}
\end{algorithm}


\end{document}